\documentclass[final]{article}

\usepackage[a4paper]{geometry}
\usepackage[T1]{fontenc}        
\usepackage{lmodern}

\makeatletter
\long\def\@makecaption#1#2{%
  \vskip\abovecaptionskip
  \sbox\@tempboxa{#1: #2}%
  \ifdim \wd\@tempboxa >\hsize
    #1: #2\par
  \else
    \global \@minipagefalse
    \hb@xt@\hsize{\hfil\box\@tempboxa\hfil}%
  \fi
  \vskip\belowcaptionskip}
\makeatother

\usepackage{amsfonts, amssymb, amsbsy, amsbsy}
\usepackage[fixamsmath,disallowspaces]{mathtools}
\usepackage{textcomp}
\usepackage{epsfig}
\usepackage{graphicx}
\usepackage{colordvi}
\usepackage{circuitikz}
\usepackage{adjustbox}
\usepackage{booktabs}         
\usepackage{array}            
\usepackage{paralist}         
\usepackage{verbatim}         
\usepackage{microtype}        
\usepackage{hyperref}
\usepackage{cleveref}         
\usepackage{tikz}
\usepackage{amsthm}
\usepackage{xspace}
\usepackage{nicefrac}
\usepackage[font=footnotesize]{caption}
\usepackage{subcaption}
\usepackage[notref,notcite]{showkeys}	
\usepackage{multirow}

\usepackage{authblk}

\usepackage{ifthen}
\usepackage{siunitx}
\usepackage{pgfkeys}
\usepackage{pgfplotstable}
\pgfplotsset{compat=1.13}
\usepackage{pgfplots}
\usepackage{tikz}
\usepgfplotslibrary{groupplots}
\usepgfplotslibrary{units}

\usepgfplotslibrary{external}
\tikzset{external/only named=true}

\edef\histoArrowHeight{9pt}

\NewDocumentCommand\HistoGetMean{mm}{
  \pgfkeys{/pgf/fpu=false}
  \pgfplotstablegetrowsof{#1}
  \pgfmathsetmacro{\LastColumn}{\pgfplotsretval - 1}
  \pgfkeys{/pgf/fpu=true}

  \pgfplotstableset{
    create on use/new/.style={
    create col/expr={\pgfmathaccuma + \thisrow{#2}}},
  }
  \pgfplotstablegetelem{\LastColumn}{new}\of{#1}
  \pgfmathsetmacro{\Sum}{\pgfplotsretval}
  \pgfmathsetmacro{\Mean}{\Sum / (\LastColumn + 1)}
}

\pgfkeys{
 /histogram/.is family, /histogram,
 default/.style =
  {%
    width = .80\linewidth, height = .65\linewidth, mean = NOTDEFINED, rescale=1,
    xmax=6e15, bins=35, hbin=70%
  },
 width/.estore in = \histoWidth,
 height/.estore in = \histoHeight,
 mean/.estore in = \histoMean,
 rescale/.estore in = \histoRescale,
 xmax/.estore in = \histoXmax,
 bins/.estore in = \histoBins,
 hbin/.estore in = \histoHBin
}

\pgfkeys{
 /histogramwawm/.is family, /histogramwawm,
 default/.style =
  {%
    width = .80\linewidth, height = .65\linewidth, rescale=1,
    xmax=6e15, bins=35, hbin=70%
  },
 width/.estore in = \histoWidth,
 height/.estore in = \histoHeight,
 rescale/.estore in = \histoRescale,
 xmax/.estore in = \histoXmax,
 bins/.estore in = \histoBins,
 hbin/.estore in = \histoHBin
}

\pgfkeys{
 /histogramwawmp/.is family, /histogramwawmp,
 default/.style =
  {%
    width = .80\linewidth, height = .65\linewidth, npoints = 1,
    hbin=70%
  },
 width/.estore in = \histoWidth,
 height/.estore in = \histoHeight,
 topmean/.estore in = \histoTopMean,
 bottommean/.estore in = \histoBottomMean,
 npoints/.estore in = \histoNPoints,
 hbin/.estore in = \histoHBin,
}

\NewDocumentCommand\DrawHistogram{oooom}{%

  \IfNoValueTF{#1}%
    {\pgfkeys{/histogram, default, }}%
    {\pgfkeys{/histogram, default, #1}}

  \edef\tableread{%
    \noexpand\pgfplotstableread{#5}{\noexpand\hdata}%
    \noexpand\pgfplotstablecreatecol[%
        create col/expr={\noexpand\thisrow{errors} / \histoRescale}%
    ]{rescalederr}\noexpand\hdata
  }

  \pgfkeys{/pgf/fpu=true}
  \tableread

  \pgfplotstablegetrowsof{\hdata}
  \edef\NTests{\pgfplotsretval}

  \pgfmathsetmacro\maxbin{\NTests * (\histoHBin / 100)}

  \pgfmathsetmacro\arrowheight{\maxbin / 10}

  \ifthenelse{%
    \equal{\histoMean}{NOTDEFINED}%
  }{%
    \HistoGetMean{\hdata}{rescalederr}%
  }{%
    \pgfmathsetmacro{\Mean}{\histoMean}%
  }
  \pgfkeys{/pgf/fpu=false}

  \begin{axis}[
      height=\histoHeight,
      width=\histoWidth,
      ymin=0, ymax=\maxbin,
      xmin=0, xmax=\histoXmax,
      ybar,
      axis x line=bottom,
      axis y line=left,
      title = {},
      xticklabel={%
        \pgfkeys{/pgf/fpu=true}%
        \pgfmathparse{\histoXmax}
        \pgfmathfloattofixed{\pgfmathresult}
        \pgfkeys{/pgf/fpu=false}
        \let\r=\pgfmathresult
        \pgfmathparse{abs(\tick - \r) < 0.001 ? int(1) : int(0)}%
        \ifthenelse{\pgfmathresult=1}{%
          \pgfmathparse{\tick}\pgfmathprintnumber{\pgfmathresult} $\! \! \times C_0$%
        }{\tick}
      },
      yticklabel={%
        \pgfmathparse{(\tick/\NTests)*100}\pgfmathprintnumber{\pgfmathresult}\%
      },
      yticklabel style={
        /pgf/number format/.cd,
        fixed, precision=0,
        /tikz/.cd
      },
      minor xtick={},
      grid=minor
  ]

  \pgfkeys{/pgf/fpu=true}
  \addplot+[hist={bins=\histoBins}, style={draw=black,fill=gray}] table[y = rescalederr]{\hdata};

  \draw[orange, ultra thick, dashed]
            (axis cs:\Mean, \pgfkeysvalueof{/pgfplots/ymin}) --
            (axis cs:\Mean, \pgfkeysvalueof{/pgfplots/ymax});

  \IfNoValueTF{#2}{}{%
    \coordinate (bluearrowstart) at (axis cs:#2, -\arrowheight);
    \coordinate (bluearrowend) at (axis cs: #2, 0);
  }

  \IfNoValueTF{#3}{}{%
    \coordinate (redarrowstart) at (axis cs:#3, -\arrowheight);
    \coordinate (redarrowend) at (axis cs: #3, 0);
  }

  \IfNoValueTF{#4}{}{%
    \coordinate (greenarrowstart) at (axis cs:#4, -\arrowheight);
    \coordinate (greenarrowend) at (axis cs: #4, 0);
  }
  \pgfkeys{/pgf/fpu=false}

  \end{axis}

  \IfNoValueTF{#2}{}{
    \draw [blue, -stealth] (bluearrowstart)  -- (bluearrowend);
  }

  \IfNoValueTF{#3}{}{%
    \draw [red, -{stealth}{stealth}] (redarrowstart)  -- (redarrowend);
  }

  \IfNoValueTF{#4}{}{%
    \draw [olive, -{stealth}{stealth}{stealth}] (greenarrowstart)  -- (greenarrowend);
  }
}

\NewDocumentCommand\DrawHistogramWithAndWithoutMean{ooooooom}{%

  \IfNoValueTF{#1}%
    {\pgfkeys{/histogramwawm, default, }}%
    {\pgfkeys{/histogramwawm, default, #1}}

  \edef\tableread{%
    \noexpand\pgfplotstableread{#8}{\noexpand\hdata}%
    \noexpand\pgfplotstablecreatecol[%
        create col/expr={\noexpand\thisrow{errors} / \histoRescale}%
    ]{rescalederr}\noexpand\hdata
    \noexpand\pgfplotstablecreatecol[%
        create col/expr={\noexpand\thisrow{errors-no-mean} / \histoRescale}%
    ]{rescalederr-no-mean}\noexpand\hdata
  }

  \pgfkeys{/pgf/fpu=true}
  \tableread

  \pgfplotstablegetrowsof{\hdata}
  \edef\NTests{\pgfplotsretval}

  \pgfmathsetmacro\maxbin{\NTests * (\histoHBin / 100)}

   \HistoGetMean{\hdata}{rescalederr}
   \pgfmathsetmacro{\TopMean}{\Mean}

   \HistoGetMean{\hdata}{rescalederr-no-mean}
   \pgfmathsetmacro{\BottomMean}{\Mean}

  \pgfkeys{/pgf/fpu=false}

  \begin{groupplot}[
    group style={
        group size=1 by 2,
        vertical sep=22pt
      },
      height=\histoHeight,
      width=\histoWidth,
      ymin=0, ymax=\maxbin,
      xmin=0, xmax=\histoXmax,
      ybar,
      axis x line=bottom,
      axis y line=left,
      title = {},
      xticklabel={%
        \pgfkeys{/pgf/fpu=true}%
        \pgfmathparse{\histoXmax}
        \pgfmathfloattofixed{\pgfmathresult}
        \pgfkeys{/pgf/fpu=false}
        \let\r=\pgfmathresult
        \pgfmathparse{abs(\tick - \r) < 0.001 ? int(1) : int(0)}%
        \ifthenelse{\pgfmathresult=1}{%
          \pgfmathparse{\tick}\pgfmathprintnumber{\pgfmathresult} $\! \! \times C_0$%
        }{\tick}
      },
      yticklabel={%
        \pgfmathparse{(\tick/\NTests)*100}\pgfmathprintnumber{\pgfmathresult}\%
      },
      yticklabel style={
        /pgf/number format/.cd,
        fixed, precision=0,
        /tikz/.cd
      },
      scaled x ticks=true,
      every x tick label/.append style={alias=XTick,inner xsep=0pt},
      every x tick scale label/.style={at=(XTick.base east),anchor=base west}
  ]

  \nextgroupplot[axis x line=bottom, xticklabel shift=3pt]
  \pgfkeys{/pgf/fpu=true}
  \addplot[hist={bins=\histoBins}, style={draw=black,fill=gray}] table[y = rescalederr]{\hdata};

  \draw[orange, ultra thick, dashed]
            (axis cs:\TopMean, \pgfkeysvalueof{/pgfplots/ymin}) --
            (axis cs:\TopMean, \pgfkeysvalueof{/pgfplots/ymax});

  \IfNoValueTF{#2}{}{%
    \coordinate (bluearrowtop) at (axis cs:#2, 0);
  }

  \IfNoValueTF{#4}{}{%
    \coordinate (redarrowtop) at (axis cs:#4, 0);
  }

  \IfNoValueTF{#6}{}{%
    \coordinate (greenarrowtop) at (axis cs:#6, - 0);
  }

  \pgfkeys{/pgf/fpu=false}

  \nextgroupplot[y dir=reverse,  axis x line=top, xticklabel=\empty, scaled x ticks=false]
  \pgfkeys{/pgf/fpu=true}
  \addplot[hist={bins=\histoBins}, style={draw=black, fill=gray}] table[y = rescalederr-no-mean]{\hdata};

  \draw[orange, ultra thick, dashed]
         (axis cs:\BottomMean, \pgfkeysvalueof{/pgfplots/ymin}) --
         (axis cs:\BottomMean, \pgfkeysvalueof{/pgfplots/ymax});

  \IfNoValueTF{#3}{}{%
    \coordinate (bluearrowbottom) at (axis cs:#3, 0);
  }

  \IfNoValueTF{#5}{}{%
    \coordinate (redarrowbottom) at (axis cs:#5, 0);
  }

  \IfNoValueTF{#7}{}{%
    \coordinate (greenarrowbottom) at (axis cs:#7, - 0);
  }

   \pgfkeys{/pgf/fpu=false}
\end{groupplot}

  \IfNoValueTF{#2}{}{
    \draw [blue, -stealth] ($(bluearrowtop) - (0, \histoArrowHeight)$)  -- (bluearrowtop);
  }

  \IfNoValueTF{#3}{}{
    \draw [blue, -stealth] ($(bluearrowbottom) - (0, -\histoArrowHeight)$)  -- (bluearrowbottom);
  }

  \IfNoValueTF{#4}{}{%
    \draw [red, -{stealth}{stealth}] ($(redarrowtop) - (0, \histoArrowHeight)$)  -- (redarrowtop);
  }

  \IfNoValueTF{#5}{}{%
    \draw [red, -{stealth}{stealth}] ($(redarrowbottom) - (0, -\histoArrowHeight)$)  -- (redarrowbottom);
  }

  \IfNoValueTF{#6}{}{%
    \draw [olive, -{stealth}{stealth}{stealth}] ($(greenarrowtop) - (0, \histoArrowHeight)$)  -- (greenarrowtop);
  }

  \IfNoValueTF{#7}{}{%
    \draw [olive, -{stealth}{stealth}{stealth}] ($(greenarrowbottom) - (0, -\histoArrowHeight)$)  -- (greenarrowbottom);
  }
}

\NewDocumentCommand\DrawHistogramPrecomputed{oooom}{%

  \IfNoValueTF{#1}%
    {\pgfkeys{/histogramwawmp, default, }}%
    {\pgfkeys{/histogramwawmp, default, #1}}

  \edef\tableread{%
    \noexpand\pgfplotstableread[skip first n=2]{#5}{\noexpand\hdata}%
  }

  \tableread

  \pgfplotstablegetrowsof{\hdata}
  \pgfmathsetmacro{\LastColumn}{\pgfplotsretval - 1}

  \pgfkeys{/pgf/fpu=true}
  \pgfmathsetmacro\maxbin{\histoNPoints * (\histoHBin / 100)}

  \pgfplotstablegetelem{\LastColumn}{xend}\of{\hdata}
  \pgfmathsetmacro{\xmax}{\pgfplotsretval}
  \pgfkeys{/pgf/fpu=false}

  \begin{axis}[
    height=\histoHeight,
    width=\histoWidth,
    ymin=0, ymax=\maxbin,
    xmin=0, xmax=0.5,
    ybar,
    compat=1.3,
    tick label style={font=\tiny},
    label style={font=\tiny},
    axis x line=bottom,
    axis y line=left,
    title = {},
    extra x ticks={0.5},
    xticklabel={%
      \pgfmathparse{abs(\tick - 0.5) < 0.001 ? int(1) : int(0)}%
      \ifthenelse{\pgfmathresult=1}{%
        $\qquad \cdot C_0$%
      }{\tick}
    },
    scaled y ticks = false,
    yticklabel={%
      \pgfkeys{/pgf/fpu=true}%
      \pgfmathparse{(\tick / \histoNPoints) * 100}\pgfmathprintnumber{\pgfmathresult}\% %
      \pgfkeys{/pgf/fpu=false}
    },
    yticklabel style={
      /pgf/number format/.cd,
      fixed, precision=0,
      /tikz/.cd
    },
    scaled x ticks=true,
    axis x line=bottom,
  ]

  \addplot+[ybar interval, style={draw=black,fill=gray}] table[x = xstart, y = errors]{\hdata};

  \draw[orange, ultra thick, dashed]
            (axis cs:\histoTopMean, \pgfkeysvalueof{/pgfplots/ymin}) --
            (axis cs:\histoTopMean, \pgfkeysvalueof{/pgfplots/ymax});

  \IfNoValueTF{#2}{}{%
    \coordinate (bluearrowtop) at (axis cs:#2, 0);
  }

  \IfNoValueTF{#3}{}{%
    \coordinate (redarrowtop) at (axis cs:#3, 0);
  }

  \IfNoValueTF{#4}{}{%
    \coordinate (greenarrowtop) at (axis cs:#4, - 0);
  }

  \draw (0,-\histoHeight) --node[fill=black] {a} (0,0);

  \end{axis}

  \IfNoValueTF{#2}{}{
    \draw [blue, -stealth, thick] ($(bluearrowtop) - (0, \histoArrowHeight)$)  -- (bluearrowtop);
  }

  \IfNoValueTF{#3}{}{%
    \draw [red, -{stealth}{stealth}, thick] ($(redarrowtop) - (0, \histoArrowHeight)$)  -- (redarrowtop);
  }

  \IfNoValueTF{#4}{}{%
    \draw [olive, -{stealth}{stealth}{stealth}, thick] ($(greenarrowtop) - (0, \histoArrowHeight)$)  -- (greenarrowtop);
  }

  \draw [black, opacity=0] ($(greenarrowtop) - (0, 3.4)$)  -- (greenarrowtop);
  
}

\NewDocumentCommand\DrawHistogramWithAndWithoutMeanPrecomputed{ooooooom}{%

  \IfNoValueTF{#1}%
    {\pgfkeys{/histogramwawmp, default, }}%
    {\pgfkeys{/histogramwawmp, default, #1}}

  \edef\tableread{%
    \noexpand\pgfplotstableread[skip first n=2]{#8}{\noexpand\hdata}%
  }

  \tableread

  \pgfplotstablegetrowsof{\hdata}
  \pgfmathsetmacro{\LastColumn}{\pgfplotsretval - 1}

  \pgfkeys{/pgf/fpu=true}
  \pgfmathsetmacro\maxbin{\histoNPoints * (\histoHBin / 100)}

  \pgfplotstablegetelem{\LastColumn}{xend}\of{\hdata}
  \pgfmathsetmacro{\xmax}{\pgfplotsretval}
  \pgfkeys{/pgf/fpu=false}

  \begin{groupplot}[
    group style={
        group size=1 by 2,
        vertical sep=14pt
      },
      height=\histoHeight,
      width=\histoWidth,
      ymin=0, ymax=\maxbin,
      xmin=0, xmax=0.5,
      ybar,
      compat=1.3,
      tick label style={font=\tiny},
      label style={font=\tiny},
      axis x line=bottom,
      axis y line=left,
      title = {},
      scaled y ticks = false,
      extra x ticks={0.5},
      xticklabel={%
        \pgfmathparse{abs(\tick - 0.5) < 0.001 ? int(1) : int(0)}%
        \ifthenelse{\pgfmathresult=1}{%
          $\qquad \cdot C_0$%
        }{\tick}
      },
      yticklabel={%
        \pgfkeys{/pgf/fpu=true}%
        \pgfmathparse{(\tick / \histoNPoints) * 100}\pgfmathprintnumber{\pgfmathresult}\% %
        \pgfkeys{/pgf/fpu=false}
      },
      yticklabel style={
        /pgf/number format/.cd,
        fixed, precision=0,
        /tikz/.cd
      },
      scaled x ticks=true,
      every x tick label/.append style={alias=XTick,inner xsep=0pt},
      every x tick scale label/.style={at=(XTick.base east),anchor=base west}
  ]

  \nextgroupplot[axis x line=bottom, xticklabel shift=0pt]

  \addplot+[ybar interval, style={draw=black,fill=gray}] table[x = xstart, y = errors]{\hdata};

  \draw[orange, ultra thick, dashed]
            (axis cs:\histoTopMean, \pgfkeysvalueof{/pgfplots/ymin}) --
            (axis cs:\histoTopMean, \pgfkeysvalueof{/pgfplots/ymax});

  \IfNoValueTF{#2}{}{%
    \coordinate (bluearrowtop) at (axis cs:#2, 0);
  }

  \IfNoValueTF{#4}{}{%
    \coordinate (redarrowtop) at (axis cs:#4, 0);
  }

  \IfNoValueTF{#6}{}{%
    \coordinate (greenarrowtop) at (axis cs:#6, - 0);
  }

  \nextgroupplot[y dir=reverse,  axis x line=top, xticklabel=\empty, scaled x ticks=false]

  \addplot[ybar interval, style={draw=black, fill=gray}] table[x = xstart, y = errors-no-mean]{\hdata};
  \draw[orange, ultra thick, dashed]
        (axis cs:\histoBottomMean, \pgfkeysvalueof{/pgfplots/ymin}) --
        (axis cs:\histoBottomMean, \pgfkeysvalueof{/pgfplots/ymax});

  \IfNoValueTF{#3}{}{%
    \coordinate (bluearrowbottom) at (axis cs:#3, 0);
  }

  \IfNoValueTF{#5}{}{%
    \coordinate (redarrowbottom) at (axis cs:#5, 0);
  }

  \IfNoValueTF{#7}{}{%
    \coordinate (greenarrowbottom) at (axis cs:#7, - 0);
  }

  \end{groupplot}

  \IfNoValueTF{#2}{}{
    \draw [blue, -stealth, thick] ($(bluearrowtop) - (0, \histoArrowHeight)$)  -- (bluearrowtop);
  }

  \IfNoValueTF{#3}{}{
    \draw [blue, -stealth, thick] ($(bluearrowbottom) - (0, -\histoArrowHeight)$)  -- (bluearrowbottom);
  }

  \IfNoValueTF{#4}{}{%
    \draw [red, -{stealth}{stealth}, thick] ($(redarrowtop) - (0, \histoArrowHeight)$)  -- (redarrowtop);
  }

  \IfNoValueTF{#5}{}{%
    \draw [red, -{stealth}{stealth}, thick] ($(redarrowbottom) - (0, -\histoArrowHeight)$)  -- (redarrowbottom);
  }

  \IfNoValueTF{#6}{}{%
    \draw [olive, -{stealth}{stealth}{stealth}, thick] ($(greenarrowtop) - (0, \histoArrowHeight)$)  -- (greenarrowtop);
  }

  \IfNoValueTF{#7}{}{%
    \draw [olive, -{stealth}{stealth}{stealth}, thick] ($(greenarrowbottom) - (0, -\histoArrowHeight)$)  -- (greenarrowbottom);
  }
}

\edef\sampleDomainlen{0.4}

\edef\sampleDomainstart{0.0}

\NewDocumentCommand\SampleGetMean{O{1}mm}{
  \pgfkeys{/pgf/fpu=false}
  \pgfplotstablegetrowsof{#2}
  \pgfmathsetmacro{\LastColumn}{\pgfplotsretval - 1}
  \pgfkeys{/pgf/fpu=true}

  \pgfplotstableset{
    create on use/new/.style={
    create col/expr={\pgfmathaccuma + (\thisrow{#3} / #1)}},
  }
  \pgfplotstablegetelem{\LastColumn}{new}\of{#2}
  \pgfmathsetmacro{\Sum}{\pgfplotsretval}
  \pgfmathsetmacro{\Mean}{\Sum / (\LastColumn + 1)}
  \pgfkeys{/pgf/fpu=false}
}

\pgfkeys{
 /sampleplot/.is family, /sampleplot,
 default/.style =
  {%
    width = .80\linewidth, height = .30\linewidth, color = red,%
    ymin = 0, ymax = 1.5e16, colsep = space, 
  },
 width/.estore in = \sampleplotWidth,
 height/.estore in = \sampleplotHeight,
 color/.estore in = \sampleplotColor,
 ymin/.estore in = \sampleplotYMin,
 ymax/.estore in = \sampleplotYMax,
 colsep/.estore in = \sampleplotColSep
}

\NewDocumentCommand\DrawSample{ommm}{%
  \IfNoValueTF{#1}%
    {\pgfkeys{/sampleplot, default, }}%
    {\pgfkeys{/sampleplot, default, #1}}

  \edef\tableread{%
    \noexpand\pgfplotstableread[%
      skip first n=2, col sep=\sampleplotColSep%
    ]{#2}{\noexpand\sampledata}%
  }

  \pgfkeys{/pgf/fpu=true}
  \tableread
  \pgfkeys{/pgf/fpu=false}

  \pgfplotstablegetrowsof{\sampledata}
  \edef\NPoints{\pgfplotsretval}

  \begin{axis}[
      width=\sampleplotWidth,
      height=\sampleplotHeight,
      axis x line=bottom,
      axis y line=left,
      ylabel = {doping}, y SI prefix=centi,  y unit={m^{-3}},%
      xlabel = {}, x SI prefix=milli, x unit={m}
    ]

    \addplot [black, dashed, opacity=1, ultra thick]%
      table[
          x expr = {\coordindex * \sampleDomainlen / \NPoints + \sampleDomainstart},
          y expr = {\thisrow{#3} / 1e6}
      ]{\sampledata};

    \addplot [\sampleplotColor, opacity=1, very thick]%
      table[
          x expr = {\coordindex * \sampleDomainlen / \NPoints + \sampleDomainstart},
          y expr = {\thisrow{#4} / 1e6}
      ]{\sampledata};

  \end{axis}
}

\NewDocumentCommand\DrawSampleWithoutMean{sommm}{%
  \edef\tableread{\noexpand\pgfplotstableread[skip first n=2]{#3}{\noexpand\sampledata}}

  \tableread

  \IfNoValueTF{#1}%
    {\pgfkeys{/sampleplot, default, }}%
    {\pgfkeys{/sampleplot, default, #2}}

  \pgfplotstablegetrowsof{\sampledata}
  \edef\NPoints{\pgfplotsretval}

  \SampleGetMean[1e6]{\sampledata}{#4}
  \pgfkeys{/pgf/fpu=true}
  \pgfmathsetmacro{\DMean}{\Mean}
  \pgfkeys{/pgf/fpu=false}

  \SampleGetMean[1e6]{\sampledata}{#5}
  \pgfkeys{/pgf/fpu=true}
  \pgfmathsetmacro{\PMean}{\Mean}
  \pgfkeys{/pgf/fpu=false}

  \begin{axis}[
      width=\sampleplotWidth,
      height=\sampleplotHeight,
      axis x line=bottom,
      axis y line=left,
      xlabel = {}, x SI prefix=milli, x unit={m}
    ]

   \addplot [black, dashed, ultra thick]%
      table[
          x expr = {\coordindex * \sampleDomainlen / \NPoints + \sampleDomainstart},
          y expr = {\thisrow{#4} / 1e6 - \DMean}
      ]{\sampledata};

    \addplot [\sampleplotColor, very thick]%
      table[
          x expr = {\coordindex * \sampleDomainlen / \NPoints + \sampleDomainstart},
          y expr = {\thisrow{#5} / 1e6 - \PMean}
      ]{\sampledata};

  \end{axis}
}

\NewDocumentCommand\DrawSampleNoLabel{ommm}{%
  \IfNoValueTF{#1}%
    {\pgfkeys{/sampleplot, default, }}%
    {\pgfkeys{/sampleplot, default, #1}}

  \edef\tableread{%
    \noexpand\pgfplotstableread[%
      skip first n=2, col sep=\sampleplotColSep%
    ]{#2}{\noexpand\sampledata}%
  }

  \pgfkeys{/pgf/fpu=true}
  \tableread
  \pgfkeys{/pgf/fpu=false}

  \pgfplotstablegetrowsof{\sampledata}
  \edef\NPoints{\pgfplotsretval}

  \begin{axis}[
      width=\sampleplotWidth,
      height=\sampleplotHeight,
      axis x line=bottom,
      axis y line=left,
      xlabel = {}, x SI prefix=milli, x unit={m}
    ]

    \addplot [black, dashed, opacity=1, ultra thick]%
      table[
          x expr = {\coordindex * \sampleDomainlen / \NPoints + \sampleDomainstart},
          y expr = {\thisrow{#3} / 1e6}
      ]{\sampledata};

    \addplot [\sampleplotColor, opacity=1, very thick]%
      table[
          x expr = {\coordindex * \sampleDomainlen / \NPoints + \sampleDomainstart},
          y expr = {\thisrow{#4} / 1e6}
      ]{\sampledata};

  \end{axis}
}

\NewDocumentCommand\DrawSampleWithAndWithoutMean{sommm}{%
  \edef\tableread{\noexpand\pgfplotstableread[skip first n=2]{#3}{\noexpand\sampledata}}

  \tableread

  \IfNoValueTF{#1}%
    {\pgfkeys{/sampleplot, default, }}%
    {\pgfkeys{/sampleplot, default, #2}}

  \pgfplotstablegetrowsof{\sampledata}
  \edef\NPoints{\pgfplotsretval}

  \SampleGetMean[1e6]{\sampledata}{#4}
  \pgfkeys{/pgf/fpu=true}
  \pgfmathsetmacro{\DMean}{\Mean}
  \pgfkeys{/pgf/fpu=false}

  \SampleGetMean[1e6]{\sampledata}{#5}
  \pgfkeys{/pgf/fpu=true}
  \pgfmathsetmacro{\PMean}{\Mean}
  \pgfkeys{/pgf/fpu=false}

  \IfBooleanTF#1{
    \pgfmathsetmacro{\groupplotWidth}{\sampleplotWidth}
    \pgfmathsetmacro{\groupplotHeight}{\sampleplotHeight / 2}
    \edef\groupshape{1 by 2}
  }{
    \pgfmathsetmacro{\groupplotWidth}{\sampleplotWidth / 2}
    \pgfmathsetmacro{\groupplotHeight}{\sampleplotHeight}
    \edef\groupshape{2 by 1}
  }

  \begin{groupplot}[
    group style={
        group size=\groupshape,
        horizontal sep=42pt
      },
      width=\groupplotWidth,
      height=\groupplotHeight,
      axis x line=bottom,
      axis y line=left,
      xlabel = {}, x SI prefix=milli, x unit={m}
    ]

    \nextgroupplot[ylabel = {doping}, y SI prefix=centi,  y unit={m^{-3}}]

    \addplot [black, dashed, ultra thick]%
      table[
          x expr = {\coordindex * \sampleDomainlen / \NPoints + \sampleDomainstart},
          y expr = {\thisrow{#4} / 1e6}
      ]{\sampledata};

    \addplot [\sampleplotColor, very thick]%
      table[
          x expr = {\coordindex * \sampleDomainlen / \NPoints + \sampleDomainstart},
          y expr = {\thisrow{#5} / 1e6}
      ]{\sampledata};

  \nextgroupplot[]

   \addplot [black, dashed, ultra thick]%
      table[
          x expr = {\coordindex * \sampleDomainlen / \NPoints + \sampleDomainstart},
          y expr = {\thisrow{#4} / 1e6 - \DMean}
      ]{\sampledata};

    \addplot [\sampleplotColor, very thick]%
      table[
          x expr = {\coordindex * \sampleDomainlen / \NPoints + \sampleDomainstart},
          y expr = {\thisrow{#5} / 1e6 - \PMean}
      ]{\sampledata};

  \end{groupplot}
}

\pgfkeys{
 /resnetScatter/.is family, /resnetScatter,
 default/.style =
  {%
    width = .80\linewidth, height = .65\linewidth, thickness=very thick%
  },
 width/.estore in = \resnetScatterWidth,
 height/.estore in = \resnetScatterHeight,
 thickness/.estore in = \resnetScatterThick
}

\NewDocumentCommand\DrawResnetScatterPlot{om}{%
  \IfNoValueTF{#1}%
    {\pgfkeys{/resnetScatter, default, }}%
    {\pgfkeys{/resnetScatter, default, #1}}

  \edef\tableread{%
    \noexpand\pgfplotstableread{#2}{\noexpand\sdata}%
    \noexpand\pgfplotstablecreatecol[%
        create col/expr={\noexpand\thisrow{linfty_error} / 1e22}%
    ]{rescalederr}\noexpand\sdata
  }

  \tableread

  \begin{axis}[
      height=\resnetScatterHeight,
      width=\resnetScatterWidth,
      axis x line=bottom,
      axis y line=left,
      ymode=log,
      xmode=log,
      title = {},
      xlabel={learning rate},
      ylabel={$\ell^\infty$ average error ($C_0$)},
      scaled x ticks=true,
  ]

    \addplot+[only marks, mark=diamond, color=blue, x filter/.code={%
      \pgfplotstablegetelem{\coordindex}{config}\of{\sdata}
      \ifnum\pgfplotsretval=2579{}\else \def\pgfmathresult{}\fi
    }] table[x = lr, y = rescalederr]{\sdata};

  \addplot+[only marks, mark=square, color=green, x filter/.code={%
      \pgfplotstablegetelem{\coordindex}{config}\of{\sdata}
      \ifnum\pgfplotsretval=3329{}\else \def\pgfmathresult{}\fi
    }] table[x = lr, y = rescalederr]{\sdata};

  \addplot+[only marks, mark=*, every mark/.append style={solid, fill=red}, color=red, x filter/.code={%
      \pgfplotstablegetelem{\coordindex}{config}\of{\sdata}
      \ifnum\pgfplotsretval=5168{}\else \def\pgfmathresult{}\fi
  }] table[x = lr, y = rescalederr]{\sdata};

  \addplot+[only marks, mark=triangle, color=gray, x filter/.code={%
      \pgfplotstablegetelem{\coordindex}{config}\of{\sdata}
      \ifnum\pgfplotsretval=6446{}\else \def\pgfmathresult{}\fi
    }] table[x = lr, y = rescalederr]{\sdata};

  \end{axis}
}

\definecolor{clr1}{RGB}{200,200,200}
\definecolor{clr2}{RGB}{160,160,160}
\definecolor{clr3}{RGB}{120,120,120}
\definecolor{clr4}{RGB}{80,80,80}
\definecolor{clr5}{RGB}{40,40,40}
\definecolor{clr6}{RGB}{0,0,0}

\NewDocumentCommand\DrawResnetTrainingVsBatches{om}{%
  \IfNoValueTF{#1}%
    {\pgfkeys{/resnetScatter, default, }}%
    {\pgfkeys{/resnetScatter, default, #1}}

  \pgfplotstableread{#2}{\sdata}

  \begin{axis}[
      height=\resnetScatterHeight,
      width=\resnetScatterWidth,
      axis x line=bottom,
      axis y line=left,
      title = {Batch size},
      xlabel={$\ell^\infty$ average error ($C_0$)},
  ]
  \foreach \bs/\op in {32/1, 64/2, 128/3, 256/4, 384/5, 512/6}{
    \edef\plotbatchdensity{%
        \noexpand\addplot+[mark=none, thick, solid, color=clr\op ] %
            table[x = xval, y = bsize\bs]{\noexpand\sdata};
    }
    \plotbatchdensity
  }
  \end{axis}
}

\NewDocumentCommand\DrawResnetTrainingVsWeightDecay{om}{%
  \IfNoValueTF{#1}%
    {\pgfkeys{/resnetScatter, default, }}%
    {\pgfkeys{/resnetScatter, default, #1}}

  \pgfplotstableread{#2}{\sdata}

  \begin{axis}[
      height=\resnetScatterHeight,
      width=\resnetScatterWidth,
      axis x line=bottom,
      axis y line=left,
      title = {Weight decay},
      xlabel={$\ell^\infty$ average error ($C_0$)},
  ]
  \foreach \wd/\col in {0.0/1, 0.0001/2, 0.001/3, 0.01/4, 0.1/5}{
    \edef\plotwdecaydensity{%
        \noexpand\addplot+[mark=none, thick, solid, color=clr\col]%
            table[x = xval, y = wdecay\wd]{\noexpand\sdata};
    }
    \plotwdecaydensity
  }
  \end{axis}
}

\NewDocumentCommand\DrawResnetTrainingVsGradientClipping{om}{%
  \IfNoValueTF{#1}%
    {\pgfkeys{/resnetScatter, default, }}%
    {\pgfkeys{/resnetScatter, default, #1}}

  \pgfplotstableread{#2}{\sdata}

  \begin{axis}[
      height=\resnetScatterHeight,
      width=\resnetScatterWidth,
      axis x line=bottom,
      axis y line=left,
      title = {Gradient Clipping},
      xlabel={$\ell^\infty$ average error ($C_0$)},
  ]
  \foreach \gc/\gcl in {0.01/0.01, 0.1/0.1, 1.0/1}{
    \edef\plotwdecaydensity{%
        \noexpand\addplot+[mark=none, \resnetScatterThick]%
            table[x = xval, y = gclip\gc]{\noexpand\sdata};
        \noexpand\addlegendentry{\gcl}
    }
    \plotwdecaydensity
  }
  \end{axis}
}

\makeatletter
\DeclareRobustCommand\onedot{\futurelet\@let@token\@onedot}
\def\@onedot{\ifx\@let@token.\else.\null\fi\xspace}
\makeatother

\newcommand{\Def}{:=}
\NewDocumentCommand{\nrml}{m}{\ensuremath{\widehat{#1}}}

\NewDocumentCommand{\frq}{}{\ensuremath{N_b}}
\NewDocumentCommand{\dopspace}{}{\ensuremath{\DopingSpace _{\frq}}}

\NewDocumentCommand{\ilength}{}{\ensuremath{k}}

\NewDocumentCommand{\Inverse}{}{\ensuremath{F}}
\NewDocumentCommand{\GlobalInverse}{}{\ensuremath{\widetilde{\Inverse}}}
\NewDocumentCommand{\Forward}{}{\ensuremath{U}}
\NewDocumentCommand{\SignalSpace}{}{\ensuremath{\mathcal{U}}}
\NewDocumentCommand{\Restriction}{}{\ensuremath{\text{TF}}}
\NewDocumentCommand{\SignalSpaceLPS}{}{\ensuremath{\mathcal{U}_{\Restriction}}}
\NewDocumentCommand{\DopingSpace}{}{\ensuremath{\mathcal{C}}}
\NewDocumentCommand{\DopingSpaceLPS}{}{\ensuremath{\mathcal{C}_{\Restriction}}}
\NewDocumentCommand{\DopingSpaceLPSuSigma}{}{\ensuremath{\mathcal{C}_{\Restriction,u,\Sigma}}}

\NewDocumentCommand{\totalResNetConfigs}{}{{7776}\xspace}
\NewDocumentCommand{\BlkFCB}{}{FixedChannel}
\NewDocumentCommand{\BlkBasic}{}{Basic}
\NewDocumentCommand{\TRUE}{}{\texttt{True}}

\newcommand{\be}{\begin{equation}}
\newcommand{\ee}{\end{equation}}
\newcommand{\bea}{\begin{eqnarray}}
\newcommand{\eea}{\end{eqnarray}}
\newcommand{\beaa}{\begin{eqnarray*}}
\newcommand{\eeaa}{\end{eqnarray*}}

\newcommand{\n}{{n}}
\newcommand{\p}{{p}}

\newcommand{\cc}{{c}}
\newcommand{\vv}{{v}}

\DeclareMathOperator*{\argmin}{arg\,min}

\definecolor{myblack}{rgb}{0.2,0.15,0.17}
\definecolor{myred}{rgb}{0.796, 0.157, 0}
\definecolor{mygreen}{rgb}{.57, .70, .28}
\definecolor{myblue}{rgb}{0.0, .44, 1.0}

\theoremstyle{definition}
\theoremstyle{plain}
\newtheorem{theorem}{Theorem}

\newtheorem{remark}[theorem]{Remark}

\graphicspath{{pics/}}

\ifpdf
\hypersetup{ pdftitle={Data-driven solutions of ill-posed inverse problems arising from doping reconstruction in semiconductors} }
\fi

\title{Data-driven solutions of ill-posed inverse problems arising from doping reconstruction in semiconductors}

\author[1]{Stefano Piani}
\author[3]{Patricio Farrell}
\author[4]{Wenyu Lei}
\author[2]{Nella Rotundo}
\author[1]{Luca Heltai}
\affil[1]{SISSA, Via Bonomea 265, 34136 Trieste, Italy}
\affil[2]{University of Florence, Viale Morgagni 67/A 50134 Florence, Italy}
\affil[3]{Weierstrass Institute (WIAS), Mohrenstr. 39, 10117 Berlin, Germany}
\affil[4]{University of Electronic Science and Technology of China, No.2006, Xiyuan Ave, West Hi-Tech Zone,
	611731 Chengdu, China }

\date{}    

\begin{document}

\maketitle

\begin{abstract}
The non-destructive estimation of doping concentrations in semiconductor devices
is of paramount importance for many applications ranging from crystal growth, the recent redefinition of the 1kg to
defect and inhomogeneity detection. 
A number of technologies (such as LBIC, EBIC and LPS) 
have been developed which allow the detection of doping
variations via photovoltaic effects. The idea is to illuminate the sample at 
several positions, and detect the resulting voltage drop or current at the contacts. 
We model a general class of such photovoltaic technologies 
by ill-posed global and local inverse problems based on a drift-diffusion system
which describes charge transport in a self-consistent
electrical field. The doping profile is included as a parametric field.
To numerically solve a physically relevant local inverse problem, we present three different data-driven approaches, based on 
least squares, multilayer perceptrons, and residual neural networks.
Our data-driven methods reconstruct the doping
profile for a given spatially varying voltage signal induced by a laser scan
along the sample's surface. The methods are trained on synthetic data
sets (pairs of discrete doping profiles and corresponding photovoltage signals at different
illumination positions) which are generated by efficient physics-preserving finite volume solutions
of the forward problem.  While the linear least square method yields an average absolute
$\ell^\infty$ error around $10\%$, the nonlinear networks roughly halve this error to $5\%$,
respectively.
Finally, we optimize the relevant hyperparameters and test the
robustness of our approach with respect to noise.

\end{abstract}




\section{Introduction}

Noninvasively estimating doping inhomogeneities in semiconductors is relevant for many industrial applications, ranging from controlling the semiconductor crystal purity during and after growth, the recent redefinition of the 1kg to detecting defects in the final semiconductor devices such as solar cells.
Doping variations lead to local electrical fields.
Experimentalists may exploit this mechanism to identify inhomogeneities, variations, and defects in the
doping profile by systematically generating electron-hole pairs via some form of electromagnetic radiation.
Due to the local fields generated by the doping inhomogeneities, the charge carrier tend to redistribute in the region
surrounding the excitation to minimize the energy. The remaining charge carriers will flow through the external circuit, and induce a current which
may be measured. By scanning the semiconductor sample with the electromagnetic source at different positions, 
one can eventually visualize the 
distribution of electrically active charge-separating defects and variations in the doping profile of the sample along the scan locations.

Several different technologies use the described photovoltage mechanism to analyze doping inhomogeneities. They are classified according to either the type of
excitation used to induce the local generation of electron-hole pairs,
or the contact placement to measure the generated current, 
or how the collected signal is related to the doping variation. 
Electron Beam Induced Current (EBIC)~\cite{Wittry1967}, Laser Beam Induced
Current (LBIC)~\cite{Bajaj1987}, scanning photovoltage
(SPV)~\cite{Jastrzebski1982}, and Lateral Photovoltage Scanning
(LPS)~\cite{Luedge1997, Kayser2020b} are some of such photovoltaic technologies, using as electromagnetic source either localized electron beams (EBIC) or laser beams (LBIC, SPV, and LPS).
A typical outcome of these techniques is an image where the intensity of each pixel is proportional to the total current signal induced by a beam shone on the pixel location.


The LPS method, proposed in~\cite{Luedge1997,Farrell2021} and schematically shown in
\Cref{fig:scheme} right panel, is especially useful in the context of crystal growth, where it
is virtually impossible to predict the quality (i.e., the symmetry) of a
semiconductor crystal \textit{during} its growth in a furnace. During the growth
process, thermal fluctuations near the the solid-liquid interface introduce
local fluctuations (or striations) in the doping profile. LPS detects such
doping inhomogeneities non-invasively at wafer-scale and room temperature. It is
especially suitable for low doping concentrations (\SI{E12} {cm^{-3}} ~to
\SI{E16}{cm^{-3}}).
%

Mathematically speaking, all of the discussed technologies result into inverse problems. 
The forward problem assumes 
that we know the doping profile and a set of laser spot positions. We then want to know 
the corresponding (laser spot dependent) photovoltage signals
at the contacts. The inverse problem, on the other hand, 
assumes we have measured photovoltage signals at the contacts for a laser scan across the sample. Then we want to reconstruct the doping profile at least in the probing region but ideally in the whole domain.
We model a general class of photovoltaic technologies introduced above 
by ill-posed global and local inverse problems based on a drift-diffusion system
which describes charge transport in a self-consistent
electrical field, the 
so-called van Roosbroeck system.
The doping profile is included as a parametric field.
We point out that a key difficulty in the 
mathematical modeling and numerical simulation is 
that the probing area where the laser or electron beam operates
is significantly smaller (even up to one spatial dimension) 
than the region in which we would like
to recover the doping profile. This means that the solution of such incomplete inverse
problems cannot be unique and we have to look for appropriate minimizers or rely on some
physically meaningful assumptions on the doping regarding periodicity or variation in 
only one spatial dimension.

Ill-posed inverse (PDE) problems have been studied for a long time from an analytical and a numerical point of view. We refer to the general overviews \cite{vogel2002computational,kaipio2006statistical}
and the references therein. Burger, Markowich and others analyzed inverse semiconductor problems \cite{Leitao2006, Burger2004}. Since the analysis is challenging, they often consider linearized/unipolar settings, no recombination/generation, only small external biases, linear diffusion as well as standard Dirichlet-Neumann boundary conditions. In \cite{Burger2004}, they discuss identifiability of the doping profile from capacitance, reduced current-voltage or laser beam induced current (LBIC) data. 
As for the doping reconstruction, previous numerical methods relied on optimization \cite{Peschka2018}, the level set \cite{Leitao2006} or the Landweber-Kaczmarz methods \cite{Burger2004}. Especially, the latter proved to be unstable and costly. For this reason we propose in this paper data-driven approaches to solve the local inverse photovoltage PDE problem.

In particular, we focus on three data-driven approaches: First, noting that 
under special conditions the photovoltage signal is related to the doping profile by linear operations, namely differentiation and convolution, we try a classical least squares approach. Second, allowing also a nonlinear relationship between doping profile and photovoltage signal, we train multilayer perceptrons. Finally, we 
adjust more advanced ResNets to our specific setup to solve the inverse photovoltage problem.
Any data-driven technique heavily relies on the amount and quality of the training data. 
Since, for example, growing crystals is exceptionally costly, we are not able to generate large real-life datasets. For this reason we generate physics-preserving synthetic data (measured signals and corresponding doping profiles) via a fast and efficient implementation of the forward PDE model \cite{Farrell2021} which relies on the Voronoi finite volume discretization described in \cite{Farrell2017}. 
As discussed above, in the context of LPS it has been shown that the forward model nicely encompasses 
three main physically meaningful features \cite{Farrell2021}. The flux discretization is handled by ideas of
Scharfetter and Gummel \cite{Scharfetter1969}. 
Compared to an implementation based on
commercial software \cite{kayser2018,Kayser2020}, our open-source code 
reduces the simulation time of the forward model by two orders of magnitude.
Finally, we will study the robustness with respect to noise and carefully tune
the hyperparameters. We introduce noise in the amplitude as well as in the wave
lengths and phase shifts of the doping fluctuations. While the results for noisy
data are worse than for clean data, our approach appears to be relatively robust
with respect to noise.

The literature on how deep neural networks maybe be used to solve PDEs has been
rapidly increasing in recent years. A notable numerical approach is called
physics informed neural networks (PINN); see \cite{raissi2019physics} and
references therein for more details. The key idea is to replace classical
statistical loss functions with PDE residuals. Unfortunately, in our
case  this is not directly feasible since solving the forward problem requires knowledge
of the doping profile in the entire three-dimensional domain. For the inverse
problem, however, the doping profile may only be reconstructed within a two- or even
one-dimensional subset. Other approaches include supervised deep learning
algorithms \cite{ray2018artificial,mishra} to efficiently approximate quantities
of interest for PDE solutions and \cite{lye2021multi} and reference therein to
accelerate existing PDE numerical schemes. In order to solve high-dimensional
PDEs, we refer to \cite{han2018solving,han2017deep}. In terms of recovering
parameters in governing PDEs by solving inverse problems, we refer to
\cite{chen2020physics,mishra2022estimates,lye2021iterative,nguyen2021model,sheriffdeen2019accelerating}.
In particular, \cite{chen2020physics} considers an inverse scattering problem
for nano-optics. How to recover generalized boundary conditions is studied in
\cite{mishra2022estimates}. Active learning algorithms are also proposed to
approximate associated PDE constrained optimization problems
\cite{lye2021iterative}. In \cite{nguyen2021model}, the authors provide
model-constrained deep learning approaches for inverse problems. Such strategies
are applied in \cite{sheriffdeen2019accelerating} to reconstruct the second-order 
coefficients in elliptic problems. Especially in the context of
semiconductors, it is often not clear which material parameters such as life
times or reference densities are actually correct. For this reason, the authors
in~\cite{Knapp2021} used machine learning techniques to estimate material
parameters in the context of organic semiconductors. 

The remainder of the paper is organized as follows: In \Cref{sec:basic_model},
we introduce general (forward) PDE models for general electromagnetic source
terms and doping profiles based on the van Roosbroeck system which describes
charge transport in a self-consistent electrical field. In \Cref{sec:inverse},
we present global and local inverse photovoltage models for generic
electromagnetic source terms. In \Cref{sec:methodologies} we present the three
data-driven approaches we use to solve the local inverse photovoltage problem
for a specific LPS setup. We conclude with a summary and an outlook in
\Cref{sec:summary}.
\section{Forward photovoltage model}
\label{sec:basic_model}
In this section, we first describe the drift-diffusion charge transport model and then the
circuit model which represents the volt meter, modeled by boundary condition. 
When combining both models, we are able to formulate the forward
photovoltage model which may be used to predict a measured signal (either a current or a
voltage) for given doping profile and source of electromagnetic radiation.
\subsection{The van Roosbroeck model}
\label{sec:vrs_model}

The semiconductor crystal is modeled by a bounded domain $\Omega\subset
\mathbb{R}^3$ in which two charge carriers evolve: electrons with negative
elementary charge $-q$, and holes with positive elementary charge $q$.  The
doping profile is given by the difference of donor and acceptor concentrations,
$N_{D}(\mathbf{x})-N_{A}(\mathbf{x})=:C(\mathbf{x})$, where
$\mathbf{x}=(x,y,z)^T\in \Omega$. We assume that $C$ is a bounded function, and
we call $\DopingSpace \subseteq L^{\infty}(\Omega)$ the space of all admissible
doping profiles. 

We describe the charge transport within the crystal in terms of the
electrostatic potential  $\psi(\mathbf{x})$, and quasi-Fermi potentials for electron and holes, $\varphi_{\n}(\mathbf{x})$ and  $\varphi_{\p}(\mathbf{x})$, respectively.
The current densities for electrons and holes are given by
$\boldsymbol{J}_\n(\mathbf{x})$, $\boldsymbol{J}_\p(\mathbf{x})$.
These variables shall satisfy the so-called van Roosbroeck model in which the first equation, a nonlinear Poisson equation is self-consistently coupled with
two continuity equations \cite{Markowich1986}
\begin{equation}
\begin{split}
-\nabla \cdot (\varepsilon \nabla \psi) &= q(p
-n
+ C(\mathbf{x})),
\\
-1/q \nabla\cdot \boldsymbol{J}_\n &= G(\mathbf{x}; \mathbf{x}_0) - H 
, \qquad
\boldsymbol{J}_\n = -q \mu_\n  n
\nabla \varphi_\n,
\\
1/q \nabla\cdot \boldsymbol{J}_\p &= G(\mathbf{x}; \mathbf{x}_0) - H
, \qquad
 \boldsymbol{J}_\p =   -q \mu_\p  p 
 \nabla \varphi_{\p}.
\end{split}
\label{eq:vR-model}
\end{equation}
In the van Roosbroeck model \eqref{eq:vR-model}, the permittivity of the medium is denoted by $\varepsilon$ and the mobilities of electrons and holes are respectively indicated by $\mu_n$ and $\mu_p$.
Assuming so-called Boltzmann statistics,
the relations between the quasi-Fermi potentials and the densities of electrons and holes,
$n$ and $p$ respectively, are given by
\begin{equation}
n
= N_\cc \exp \left( \frac{q(\psi - \varphi_\n)-E_{\cc}}{k_B T}\right)
\qquad \text{and} \qquad
p
= N_\vv \exp \left( \frac{ q(\varphi_\p - \psi) + E_{\vv}}{k_B T}\right).
\label{eq:dens-pot}
\end{equation}
Here, we have denoted the conduction and valence band densities of states with $N_\cc$ and $N_\vv$,
the Boltzmann constant with $k_B$ and the temperature with $T$.
Furthermore, $E_{\cc} $ and $E_{\vv} $ refer to the constant  conduction
and valence band-edge energies, respectively.

The semiconductor is considered in \textit{equilibrium} if the quasi-Fermi potentials vanish, $\varphi_n=\varphi_p=0$. In this case only the nonlinear Poisson equation in \eqref{eq:vR-model} 
is solved for the equilibrium electrostatic potential $\psi_{eq}$. The corresponding equilibrium charge carrier densities $n_{eq}$ and $p_{eq}$ satisfy $n_{eq}p_{eq} = n_i^2$, where
$n_i$ is the so-called intrinsic carrier density, defined via the relationship
$n_i^2 = N_c N_v\exp\left(-{(E_c-E_v)}/{k_B T}\right)$.

The recombination term $H$ is the sum of the direct recombination, the Auger recombination and the Shockley-Read-Hall recombination, that is respectively: 
\begin{equation*}
  H_{\si{dir}}=C_d(np-n_i^2),
\quad 
  H_{\si{Aug}}=C_{\n}n(np-n_i^2)+C_{\p}p(np-n_i^2),
\quad 
  H_{\si{SRH}}=\frac{np-n_i^2}{\tau_p(n+n_T)+\tau_n(p+p_T)}.
\end{equation*}

Different types of crystal samples (such as silicon, germanium, gallium arsenide) have different life times $\tau_n, \tau_p$, and reference densities $n_T, p_T$. 

\begin{figure}[!htb]
\begin{subfigure}{.47\textwidth}
\begin{circuitikz}
  \draw (0,0) node[ground] {}
              node[anchor=east] {$u_{ref}$}
              to[short,*-*]
        (0,1.5) node[anchor=south east] {$u_{D_1} $}
              to[photodiode,-,color=red]
        (5,1.5) node[anchor=south west] {$u_{D_2} $}
              to[short,*-]
        (5,0) to[R, l_=$R$, i>_=$i_{D} $] (0,0);

  \begin{scope}[shift={(0,-3)}]

        \draw
            (0,1) to[photodiode,o-o,color=red]
            (1,1) node[anchor=west] {$:=$};

        \draw[line width=2]
            (2.5,0) -- node[anchor=east] {$\Gamma_{D_1}$} (2.5,2);
            
        \draw
            (2.5,2) -- node[anchor=south] {$\Gamma_{N}$} (5,2);
        
        \draw
            (5,0) -- node[anchor=north] {$\Gamma_{N}$} (2.5,0);

        \draw[line width=2]
            (5,2) -- node[anchor=west] {$\Gamma_{D_2}$} (5,0);

        \node at (3,.5) {$\Omega$};
        
         \draw[line width=1.5]
            (3.75,1) -- node[anchor=south] {$\Sigma$} (4.7,1);
  \end{scope}
\end{circuitikz}
\end{subfigure}
\hfill
\begin{subfigure}{.52\textwidth}
    \includegraphics[width=.98\textwidth]{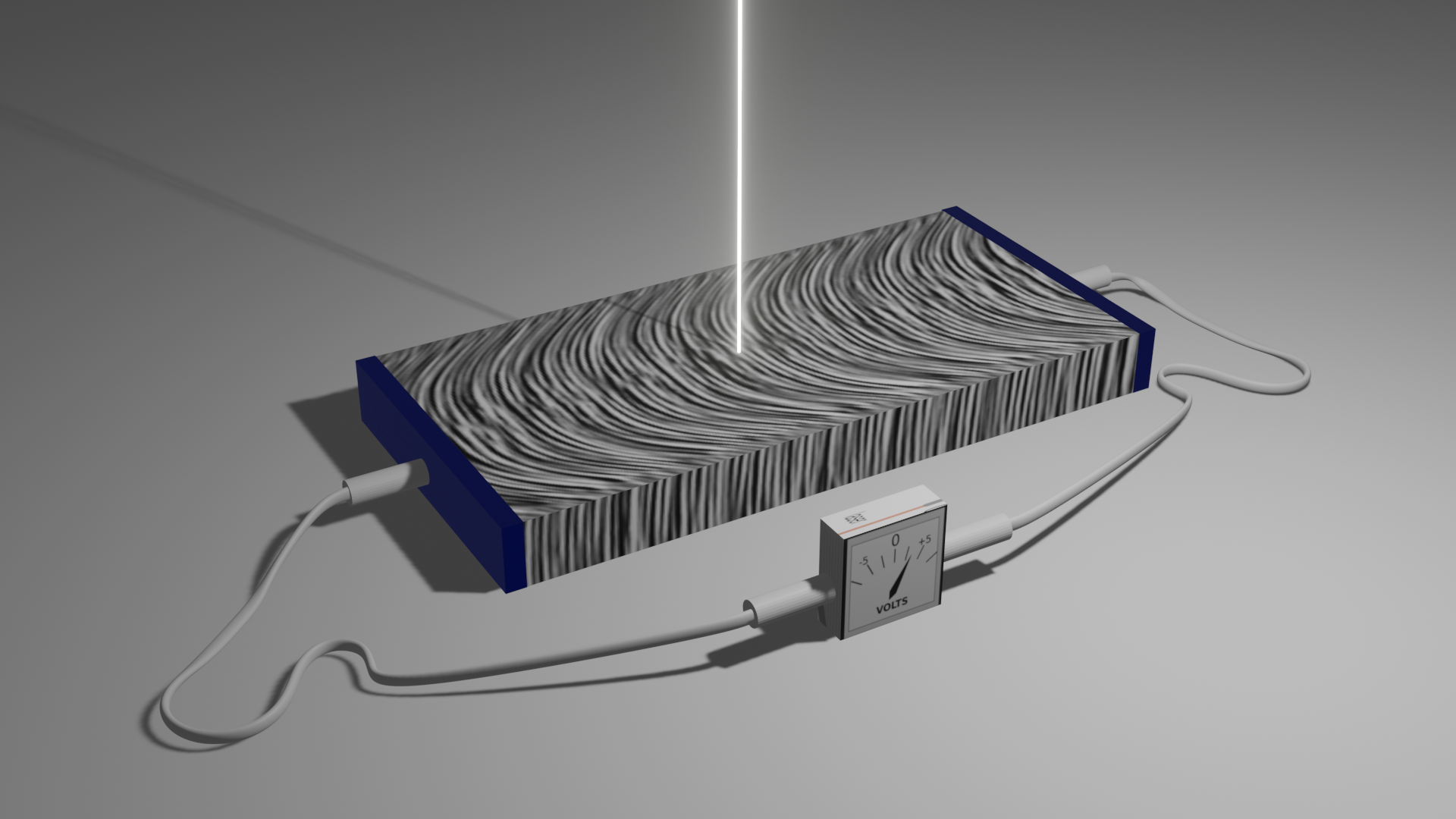}
\end{subfigure}
\caption{%
  On the top left we show the schematic of a photo-sensitive silicon crystal (in red) coupled with the voltage
  meter having resistance $R$ (rendered on the right). On the bottom left, we show the domain $\Omega$, the ohmic contacts $\Gamma_{D_1}$ and $\Gamma_{D_2}$, and the laser probing area $\Sigma$.
}
\label{fig:scheme}
\end{figure}
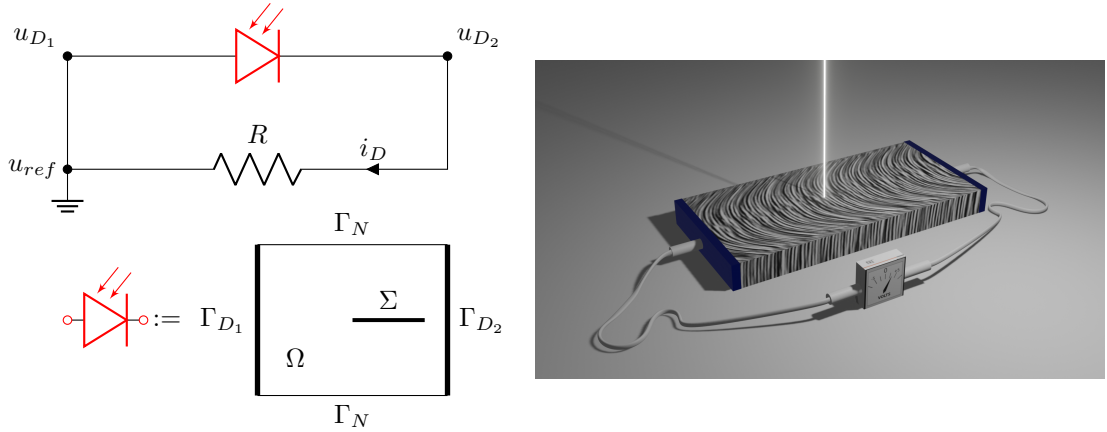

The electromagnetic source (a laser or an electron beam) is modeled by the generation term $G(\mathbf{x}; \mathbf{x}_0)$. When the laser hits the crystal at the point $\mathbf{x}_0 := (x_0, y_0, z_0)^T$, some photons are \textit{impinged} and create electron-hole pairs, resulting in a  generation rate defined as follows
\begin{equation}
G(\mathbf{x}; \mathbf x_0)=\kappa S(\mathbf{x}-\mathbf{x}_0),
\label{eq:generation-rate}
\end{equation}
where $S(\mathbf{x})$ is the shape function of the laser (normalized by $\int_{\mathbb{R}^3}
S(\mathbf{x}){d}\mathbf{x}=1$), while $\kappa$ is a constant given by
$\kappa := \frac{P\lambda_L}{h}(1-r)$. 
Here, $P$ is the laser power, $\lambda_{{L}}$ is the wave length of the laser, $h$ is the Planck constant, and $r$ is the reflectivity rate of the crystal.

\begin{figure}[!htb]
  \begin{center}
   \includegraphics
   {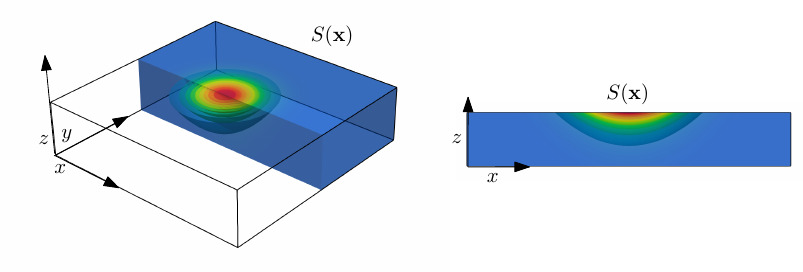}
  \end{center}
  \caption{A zoom-in of the area around the point $\mathbf{x}_0$, where the laser hits the crystal (left), and a cross section of the device at $y=0$ used in the 2D simulation (right).}
  \label{fig:laser_shape}
\end{figure}

We assume that area of influence of the electromagnetic source decays exponentially fast from the incident point $\mathbf{x}_0$. In particular, we take a laser profile function $S$ defined as
\begin{equation}
    S(\mathbf{x}) := \frac{1}{2\pi\sigma_L^2 d_A}
    \exp\left[-\frac12 \left( \frac{x}{\sigma_L}\right)^2 \right]
    \exp\left[-\frac12 \left( \frac{y}{\sigma_L}\right)^2 \right]
    \exp\left[-\frac{|z|}{d_{A}}\right].
\label{eq:shapeS}
\end{equation}
Here $\sigma_L$ is the laser spot radius, while $d_{A}$ is the penetration depth
(or the reciprocal of the absorption coefficient), which  heavily depends on the
laser wave length. \Cref{fig:laser_shape} shows a typical configuration, where
the laser beam hits the crystal on the top surface, and shows the exponential
decay of the laser shape function $S$. 

A prototypical setting for photovoltage measurements is given by a cuboid sample attached to a
voltmeter with resistance $R$, and the electromagnetic source aiming at its top surface.

%
\subsection{Boundary conditions}
The PDE system \eqref{eq:vR-model} is supplemented with Dirichlet and Neumann boundary conditions.
The boundary  $\partial\Omega$ is the union of two disjoint parts $\Gamma_N$ and $\Gamma_D$. On $\Gamma_N$,  we assign Neumann boundary conditions
\begin{equation}
\frac{\partial \psi}{\partial \boldsymbol{\nu}}=
\frac{\partial \varphi_\n}{\partial \boldsymbol{\nu}}=
\frac{\partial \varphi_\p}{\partial \boldsymbol{\nu}}=0,  
\label{eq:neumann-data}
\end{equation}
where $\partial/\partial \boldsymbol{\nu} =\boldsymbol{\nu} \cdot \nabla$ is the normal derivative along the external unit normal $\boldsymbol{\nu}$.
On $\Gamma_D$, we assign Dirichlet-type boundary conditions. This type of boundary condition models
so-called ohmic contacts \cite{Selberherr1984}. 
We suppose that there are two ohmic contacts, i.\ e. $\Gamma_D= \Gamma_{D_1} \cup \Gamma_{D_2}$. The ohmic boundary conditions can be summarized by
\begin{equation}
\psi - \psi_0= \varphi_\n = \varphi_\p = u_{D_i} - u_{\text{ref}}\quad \text{on }\quad\Gamma_{D_i}, \quad i=1,2,
\label{eq:Dirichlet-data}
\end{equation}
where $ \psi_0$ is the local electroneutral potential which one obtains by solving the Poisson equation for $\psi$ for a vanishing left-hand side.
The terms $u_{D_i}$ denote the contact voltages at the corresponding ohmic contacts; for the 
forward model their difference is a-priori unknown.
We define a reference value $u_{\text{ref}}$ of the potential and set it to
$u_{\text{ref}} = 0 = u_{D_1}$.

The total electric current $i_D $  flowing through the ohmic contact $\Gamma_{D_2}$ is defined by the surface integral
\begin{equation}
 \begin{aligned}
   i_D \colon\quad& \Omega \times \mathbb{R} \times \DopingSpace &\to &\qquad  \mathbb{R} \\
   &(\mathbf{x}_0, u_{D_2}, C)& \mapsto & \quad i_D(\mathbf{x}_0, u_{D_2}, C) :=  \int_{\Gamma_{D_2}}   \boldsymbol{\nu} \cdot (\boldsymbol{J}_\n (\mathbf{x})+\boldsymbol{J}_\p (\mathbf{x} ))  d\sigma(\mathbf{x}).
 \end{aligned}
\label{eq:current-i-th-ohmic}
\end{equation}

The current $i_D$ is a function of the contact voltage $u_{D_2}$, of the laser
spot position $\mathbf{x}_0$, and of the doping profile $C$. 

In order to close the system, we model the voltage/ampere meter as a simple
circuit having a resistance $R$. This structure is visualized in
\Cref{fig:scheme}. When the external circuit is connected to the crystal, the
boundary condition $u_{D_2}$ can no longer be chosen arbitrarily and must
balance the voltage difference at the voltmeter due to the current $i_D$, which
is equal to the product $R~i_D$. A generalized theory of this coupling can be
found in \cite{Ali2010} and in the references therein.

Let $\SignalSpace$ be the function space of all possible photovoltage signals as
a function of the laser spot position $\textbf{x}_0$ that can be measured at the contacts for
any given admissible doping profile $C\in\DopingSpace$. We define the \emph{laser-voltage} (L-V)
map as precisely the map that associates to each laser spot location
$\mathbf{x}_0$ and doping profile $C$ the corresponding photovoltage signal
$u_P$ which satisfies the second ohmic boundary condition in \eqref{eq:Dirichlet-data}
 \textit{and} balances the voltage difference across $R$, that is
\begin{equation}
  \begin{aligned}
    u_P \colon\quad& \Omega \times \DopingSpace &\to &\qquad  \mathbb{R} \\
         &(\mathbf{x}_0, C)& \mapsto & \qquad u \in \SignalSpace \text{ such that } u(\mathbf{x}_0) 
         - R~i_D( \mathbf{x}_0, u, C)= 0 \text{ and  } u(\mathbf{x}_0) = u_{D_2} .
  \end{aligned}
  \label{eq:laser-voltage}
\end{equation}
The laser-voltage map \eqref{eq:laser-voltage} corresponds to a nonlinear version of Ohm's law at the contact and constitutes an implicit equation for the Dirichlet boundary condition $u_{D_2}$
 of the van Roosbroeck system \eqref{eq:vR-model}.
An existence result for solutions of the laser-voltage map is provided
in~\cite{Busenberg1993}, with the assumption that the generation rate $G$ is
small enough. For a spatially varying doping profile $C$, the
laser-voltage map $u_P(\mathbf{x}_0, C)$ may be different for each laser spot location $\mathbf{x}_0$.

\subsection{The global forward photovoltage problem}

Although the presented forward PDE model is well-posed for low laser power
intensity and for all laser spot positions $\mathbf{x}_0 \in \Omega$, not all
positions in $\Omega$ can be illuminated by a laser. The major limiting factor
is the fact that the laser power intensity decays exponentially away from the surface of the crystal sample, leaving in fact as the
only feasible positions those that are near or on the surface of the crystal
sample.

We denote with $\Sigma \subseteq \overline{\Omega}$ the set of all admissible laser spot
positions $\mathbf{x}_0$.
Our \emph{global forward photovoltage problem} then associates to each doping profile $C$ the
function $u$ that maps each laser spot position $\mathbf{x}_0 \in \Sigma$ to the
corresponding photovoltage signal $u_P(\mathbf{x}_0, C)$, i.e., 
\begin{equation}
  \begin{aligned}
   \Forward :\quad& \DopingSpace &\to &\qquad  \SignalSpace \\
         & C & \mapsto & \qquad u: \Sigma \to \mathbb R, \qquad  u(\mathbf{x}_0) = u_P(\mathbf{x}_0, C), \qquad \mathbf{x}_0 \in \Sigma.
  \end{aligned}
  \label{eq:forward-problem}
\end{equation}


\section{Inverse photovoltage problems}
\label{sec:inverse}


From an industrial perspective, more interesting than the forward photovoltage problem is the inverse one.  How do doping inhomogeneities influence the measured voltage difference $u_P$, defined in the previous section? 
Suppose we have measured the photovoltage signal for several different laser
spot positions $\mathbf{x}_0$, how does the doping profile look like that leads
to this signal? Since the doping profile enters as a volumetric
source term defined in the whole domain in the van Roosbroeck system
\eqref{eq:vR-model}, but we can only probe part of the domain with the laser
signal, answering such a question implies solving an ill-posed inverse problem:
different global doping profiles may lead to the \textit{same} signal. Hence, we will first
formulate an idealized global inverse photovoltage problem 
and then a practically relevant local inverse photovoltage problem.

\subsection{Global inverse photovoltage problem}

Ideally, we would like to find the \textit{global doping reconstruction
operator}
\[
    \GlobalInverse := \Forward^{-1} \colon \SignalSpace \to \mathcal{P}(\DopingSpace) ,
\]
that, for a given photovoltage signal measurement $u$, returns the preimage 
$\Forward^{-1}(u)$. Here we indicate with $\mathcal{P}(\DopingSpace)$ the power set of
$\DopingSpace$, i.e., the set of all possible subsets of
$\DopingSpace$. The \textit{global inverse photovoltage problem} reads
\begin{equation}
   \GlobalInverse(u) = \DopingSpace_u, \qquad \DopingSpace_u := \{ C \in \DopingSpace \text{ such that } u_P(\cdot, C) = u \},
    \label{eq:inverse-operator} 
\end{equation}
that is, $\GlobalInverse$ is a function from $\SignalSpace$ to
$\mathcal{P}(\DopingSpace)$ (i.e., a multi-valued function of the signal $u$)
such that $C\subseteq \GlobalInverse (u_P(\cdot,C))$ for all $C \in
\DopingSpace$.

\subsection{Local inverse photovoltage problem}
\label{sec:local_inverse_problem}

In practice, however, the construction of the operator $\GlobalInverse$ (even
the construction of an approximation of $\GlobalInverse$) is nontrivial, and we
would like to simplify our problem by building a local inverse problem that is better posed.
A key point when defining a \emph{local inverse problem}, is the
fact that we cannot probe the entire domain $\Omega$
with the laser, but only the subset of all possible laser spot positions $\Sigma$ near or on the surface of $\Omega$. In the following subsections, we
will make some assumptions on the structure of  technologically feasible doping profiles. Technological feasibility depends on 
the growth process, on the technology used to inject doping in the crystal, and on the specific photovoltage technology. 
Here, we assume variation only in the $x$-direction. That is, $C(\mathbf{x})=C(x)$.
This class of doping profiles is relevant, for example, for crystal growth, where
striations in the doping profile indicate fluctuations of the temperature field during
the growth process, which are dominant in the growth direction.

We define the set of technologically feasible (TF) doping profiles and the set of corresponding signals as
\begin{equation}
  \label{eq:allowed-doping-space}
  \DopingSpaceLPS = \{ C \in  \DopingSpace \text{ such that } C(\mathbf{x}) = C(x) \}, \qquad \SignalSpaceLPS := u_P(\cdot, \DopingSpaceLPS).
\end{equation}

With these assumptions on the doping profile, it is reasonable to presume that
the photovoltage signal will be influenced more by variations of the doping
profile in the vicinity of the laser spot position, which is where most of the
charge carriers are generated, and therefore we define the \textit{local inverse
photovoltage problem} by restricting the global 
inverse problem~\eqref{eq:inverse-operator} both in terms of technologically feasible doping profiles, as
well as in terms of probing domain:
\begin{equation}
  \Inverse(u) := (\GlobalInverse(u)\cap \DopingSpaceLPS)|_{\Sigma} = (\DopingSpace_u\cap\DopingSpaceLPS)|_{\Sigma} =: \DopingSpaceLPSuSigma,
  \label{eq:inverse-operator-of-restriction} 
\end{equation}
where the restriction operator $|$ is applied to each element in the set.
The goal of the local inverse problem is to reconstruct 
doping profiles only in the probing area $\Sigma$ and not on the whole domain $\Omega$.
However, the underlying assumption is that the local doping reconstruction 
will give us information in a neighborhood of $\Sigma$, or even on all of
$\Omega$, when one makes additional assumptions (for example on doping
periodicity).
In general, the set $\DopingSpaceLPSuSigma = \Inverse(u)$ is much smaller than
the set $\DopingSpace_u:=\GlobalInverse(u)$, since it restricts the family of
admissible doping profiles to $\DopingSpaceLPS$, and discards all of the
information outside of $\Sigma$. Two doping profiles in $\DopingSpace_u$
correspond to the same element in $\DopingSpaceLPSuSigma$ as soon as
their restrictions on $\Sigma$ coincide.

In principle, it should be possible to formulate a well-posed local inverse
problem by further reducing the admissible doping profiles $\DopingSpaceLPS$ until
the set $\DopingSpaceLPSuSigma$ contains just one single element 
for
any admissible input signal $u \in \Forward(\DopingSpaceLPS)$. We do not know,
however, if this procedure is feasible, and what the necessary and
sufficient conditions on $\DopingSpaceLPS$ and on $\Sigma$ are to ensure that the
local inverse problem $\Inverse(u)$ is well-posed. 
We leave these questions for future investigations, and concentrate on numerical
approximations of the local inverse problem based on well-posed finite
dimensional approximations of $\Inverse$. In \Cref{sec:methodologies} we
propose three different strategies (of increasing complexity) to build an
approximate inverse operator $\Inverse_h$ starting from a collection of doping
profiles restricted to a discrete set of probing points $\Sigma_h$ and 
their corresponding discrete photovoltage signals.

\begin{remark}
  The numerical approximations that we construct always produce a
  unique answer for each finite sampling of a signal $u$ in $\SignalSpaceLPS$.
  We do not attempt to construct the full set $\DopingSpaceLPSuSigma$ of all
  possible doping profiles that would generate $u$, but only provide the
  sampling of a single doping profile $C$, which is, in some sense,
  a representative of the set $\DopingSpaceLPSuSigma$. The ill-posedness of the
  local inverse problem and the impossibility to recover the full doping profile
  $C$ after having solved the approximate inverse problem is one of the
  reasons why we cannot use more modern techniques (such as Physics Informed
  Neural Networks (PINN)~\cite{raissi2019physics}) that would exploit the
  residual of the forward problem to improve the construction of an approximate
  inverse operator.
\end{remark}

\section{Data-driven approximation of the inverse photovoltage problem}
\label{sec:methodologies}

Inverse problems for charge transport equations have often been tackled with
standard techniques (see, for
example,~\cite{Burger2004,Leitao2006,Peschka2018}). However, in other fields such as
image recognition or weather prediction, data-driven approaches have gained
significant momentum to solve a variety of inverse problems (see, for
example,~\cite{arridge2019solving, li2020nett, adler2017solving}). 

In order to formulate a data-driven approach for the numerical approximation of
the operator $\Inverse$, we leverage the well-posedness of the forward problem
$\Forward$ and its discrete approximation described in
\cite{Farrell2021} to formulate a discrete inverse problem
$\Inverse_h$ as a well-posed minimization problem in a finite dimensional space.

\subsection{Discrete local inverse problem}
\label{sec:discrete-local-inverse}

Let $\Sigma \subseteq \overline{\Omega}$ be a subdomain. We assume we sampled
both the photovoltage signal $u_P$ and the doping profile $C$ at
$n\in\mathbb{N}$ discrete laser spot positions which shall be contained in the
discrete mesh $\Sigma_h\subseteq \Sigma$. With the boldface symbols $\mathbf{u}$
and $\mathbf{C}$ we indicate the discrete samplings of the signal $u_P(\cdot,C)$
and the doping profile $C$ evaluated on $\Sigma_{h}$, respectively. Notice that
the generation of a single pair of signal and doping samples $(\mathbf{u},
\mathbf{C})\in \mathbb{R}^{n\times 2}$ from an admissible doping profile $C\in
\DopingSpaceLPS$ requires in fact the solution of $n$ discrete van Roosbroeck
systems~\eqref{eq:vR-model}, one for each laser spot position $\mathbf{x}_0 \in
\Sigma_{h}$.

The \textit{local discrete inverse problem} $\Inverse_h$ can be interpreted as a
(generally nonlinear) function that maps a vector of signal measurements
$\mathbf{u} \in \mathbb{R}^n$ to a vector of doping profile values $\mathbf{C}
\in \mathbb{R}^n$:
\begin{equation}
  \begin{aligned}
    \Inverse_h  \colon& \mathbb{R}^n & \to & \qquad \mathbb{R}^n \\
    & \mathbf{u} & \mapsto &\qquad \mathbf{C}.
  \end{aligned}
  \label{eq:discrete-inverse}
\end{equation}
Ideally, we would like to build $\Inverse_h$ in such a way that 
for all $u$ in
$\SignalSpaceLPS$, there exists a $C \in \DopingSpaceLPSuSigma=\Inverse(u)$
such that $\Inverse_h( \mathbf{u} := u|_{\Sigma_h}) = \mathbf{C} :=
C|_{\Sigma_h}$. However, this procedure suffers from the non-uniqueness of 
$\Inverse(u)$, which we mitigate by constructing $\Inverse_h$ through a 
well-posed minimization problem.

In practice, we build $\Inverse_h$ by minimizing the mean square error loss
\begin{equation}\label{eq:msel}
    \text{MSE}\left(\{\mathbf{u}_j, \mathbf{C}_j\}_{j=1}^{N_{\text{train}}}\right):= \frac{1}{N_{\text{train}}}
        \sum_{j=1}^{N_\text{train}} \|\Inverse_h(\mathbf u_j)-\mathbf{C}_j\|^2_{\ell^2} ,
\end{equation}
on a given \emph{training dataset} $\{\mathbf{u}_j,
\mathbf{C}_j\}_{j=1}^{N_{\text{train}}}$ consisting of a collection of
$N_{\text{train}}$ pairs of signal samples $\mathbf{u}_j$ and doping samples
$\mathbf{C}_j$. The generation process of the data is detailed in \Cref{sec:data_generation}.

The quality of the approximation of $\Inverse_h$ is evaluated by measuring the
prediction capabilities of $\Inverse_h$ on a \emph{test dataset} that is
independent of the training dataset. For the evaluation of the prediction
quality, we use the $\ell^\infty$ norm, and
define the \emph{test error} as
\begin{equation}
  \label{eq:testing-error}
    \text{TE}\left(\{\mathbf{u}_j, \mathbf{C}_j\}_{j=1}^{N_{\text{test}}}\right):= \frac{1}{N_{\text{test}}} \sum_{j=1}^{N_{\text{test}}}\|\Inverse_h(\mathbf u_j)-\mathbf C_j\|_{\ell^\infty}.
  \end{equation}
  It is worth noting that, while the definition of the $\textrm{MSE}$ function relies on the
  $\ell^2$ norm, in \cref{eq:testing-error} we are using the $\ell^\infty$ norm.
  There are several reasons for this.
  First of all, minimizing the $\ell^\infty$ error is more difficult than minimize the $\ell^2$
  error; therefore, when we build our models, we take advantage of the smoothness of the
  $\ell^2$ norm. But the $\ell^\infty$ error has several advantages: its independence from the
  size of the domain and its physical meaning of being an upper bound for the pointwise error.

  On the other hand, the $\ell^\infty$ error has several advantages: it is
  independent of the size of the domain and
  it has the physical meaning of being an upper
  bound on the error on each point. Finally, we perform some hyperparameter
  tuning on our models (i.e., given a family of models, we choose the one that performs
  better; see \Cref{ssec:res-mlp,ssec:res-ResNet}). Since we generate
  our models by minimizing the $\ell^2$ error while the tuning of the hyperparameters minimizes
  the error in $\ell^\infty$, we keep both of them under control (and therefore,
  by interpolation, we control every other error in $\ell^k$ for $k \geq 2$).
  Indeed, the training of a specific model minimizes the $\ell^2$ error without taking
  into account the effect this procedure has in the $\ell^\infty$ space; on the other hand,
  if this effect is too detrimental for the $\ell^\infty$ error, we discard that model during
  the tuning of the hyperparameters.

We develop three approaches to build $\Inverse_h$:
\begin{itemize}
  \item \emph{Least squares} (LS): we approximate $\Inverse_h(\mathbf{u})$ by the matrix vector product of an
  $n\times n$ matrix $\mathbf{A}_h$ with $\mathbf{u}$. We find the matrix by minimizing the MSEL error defined in \eqref{eq:msel}
  over all $n\times n$ matrices;

  \item \emph{Multilayer perceptron} (MLP): we approximate $\Inverse_h(\mathbf{u})$ by a
  multilayer perceptron, and use Stochastic Gradient Descent (SGD) to minimize
  the MSE;

  \item \emph{Residual neural network} (ResNet): we approximate $\Inverse_h(\mathbf{u})$
  by a residual neural network, and use SGD to minimize the MSE.
\end{itemize}

\subsection{An industrial application: the LPS setup}
\label{sec:setup}

We introduce a specific LPS setup for a silicon sample, discuss how and why we
generate the data from the forward model as well as solve the corresponding
local inverse photovoltage problem via the three data-driven approaches
introduced in \Cref{sec:discrete-local-inverse}.

We consider a silicon cuboid of the form $\Omega:=[-\nicefrac{\ell}{2},
\nicefrac{\ell}{2}]\times[-\nicefrac{w}{2},\nicefrac{w}{2}]\times[-h,0]$ with
length $\ell = \SI{3}{mm}$, width $w = \SI{0.5}{mm}$ and  height $h =
\SI{5e-5}{mm}$, ohmic contacts at $x = -\nicefrac{\ell}{2}$ (corresponding to
$\Gamma_{D_1}$) and $x = \nicefrac{\ell}{2}$ (corresponding to as
$\Gamma_{D_2}$), see \Cref{fig:scheme}. We focus on a lateral photovoltage
scanning (LPS) table-top setup, see \Cref{fig:scheme} right panel. The silicon parameters
needed for the corresponding forward model from \Cref{sec:basic_model}
can be found in \cite{Farrell2021}.
We simulate the charge transport in the two-dimensional plane where $y=0$ and reconstruct the doping profile along a line, namely the 
subdomain $\Sigma=\{(x,0,0)^T\in \Omega: x \in [0,\nicefrac{\ilength} 2]\}$ 
with $\ilength < \ell$, ensuring we are sufficiently far away from the
boundaries $x=\pm \nicefrac{\ell}{2}$ to reduce boundary effects. 
The laser scan then produces a uniform partition $0=x_1< \ldots <
x_{n}=\nicefrac{\ilength}{2}$ with mesh size $\Delta x$. The resulting laser spot positions are given by $\Sigma_{h}:=\{\mathbf
x_i\}_{i=1}^{n}=\{(x_i,0,0)^T\}_{i=1}^{n}$. In the following numerical simulations, we set
$n=1200$ and $\ilength = \SI{0.4}{mm}$.

\subsection{Data generation}
\label{sec:data_generation}

For our data-driven approaches, we need a suitable amount of data;
unfortunately, generating enough data from real life requires large budgets. Instead, we
will generate synthetic data from the forward model that we have described in
\Cref{sec:basic_model}. This data is physically meaningful in the sense that our
finite volume approximation of the LPS problem correctly incorporates physically
meaningful behavior, namely i) the signal intensity depends linearly on local
doping variations for low laser powers, ii) the signal saturates for higher
laser intensities due to a screening effect, and iii) the signal depends
logarithmically on moderate laser intensities~\cite{Farrell2021}.

As pointed out in \Cref{sec:local_inverse_problem}, we consider a family of
doping profiles that vary only along the $x$ direction, and that are constant in
both $y$ and $z$ directions.
%
To generate our synthetic datasets, we need an algorithm that produces
fluctuating doping profiles. Since within the LPS framework the doping profile
is roughly periodic, see \cite{Farrell2021}, we assume a superposition of sinusoidal
functions from which we randomly sample physically reasonable amplitudes and
wavelengths.
Therefore, 
 we define
\begin{equation}\label{eq:doping_space}
C \big(
  \mathbf{x}=(x,y, z)^T;\boldsymbol{\beta}=(C_0,\boldsymbol{\alpha},\boldsymbol{\lambda}
)^T\big)  =  C_0 \Bigg( 1 + \sum_{i = 1}^{\frq} \alpha_i \sin\frac{2\pi x}{\lambda_i} \Bigg) ,
\end{equation}
where $\boldsymbol{\beta}$ is a vector of parameters that includes a fixed average doping value $C_0=\SI{1.0E16}{cm^{-3}}$ (a typical value for silicon crystals), amplitudes $\boldsymbol\alpha:=\{\alpha_i\}\subset \{0\}\cup [0.05, 0.2]$, and wavelengths $\boldsymbol\lambda:=\{\lambda_i\}\subset [\SI{10}{\mu m}, \SI{1000}{\mu m}]$. In our simulation, we set the number of sine terms to $N_b=5$, which seems to strike a good balance between complexity and real-life situations for striations in doping profiles.

To generate an element $\{\mathbf{u}, \mathbf{C} \}$ of our dataset, we randomly
sample a parameter $\boldsymbol \beta$ with which we generate a continuous
function $C(\mathbf{x}, \boldsymbol\beta)$ and we
compute the discrete counterpart of $ u \Def U(C) $ defined in \Cref{eq:forward-problem} using the finite volume scheme developed in \cite{Kayser2021}. Finally, we restrict both $C$ and $u$ to the discrete laser spot mesh $\Sigma_h$.

In realistic physical scenarios, the doping profiles contain noises,  and they
cannot be represented exactly by \Cref{eq:doping_space}. To simulate such scenario, we perturb some of the doping profiles and generate an additional \emph{noisy}
dataset, used to test the robustness of our networks in the presence of noise. In what follows, we say that a dataset is ``noisy'' or ``clean'' if $C$ has or has not been perturbed. The algorithm used to generate noisy datasets is described
in \Cref{sec:noise}.

In the following, we summarize the datasets used to construct and test the learning models introduced in Section~\ref{ssec:res-linear}--\ref{ssec:res-ResNet}. The datasets are available on the repository~\cite{ZenodoRepo2022}.

\subsubsection*{CleanDataSET} Consisting of:\\
\indent $\bullet$ \textbf{CleanTrainingDataSET}: a dataset with
    $N_{\text{train}}=22,500$ clean samples used to train our initial model;

$\bullet$ \textbf{CleanTestDataSET}: a dataset with $N_{\text{test}}=7,500$
    clean samples. With this dataset we test the performance of
    our models trained by \textbf{CleanTrainingDataSET}. Note that testing errors generated by this dataset (cf. \eqref{eq:testing-error}) implies the quality of a trained model. This will help use tune the corresponding hyperparameters;

$\bullet$ \textbf{CleanValidationDataSET}: a dataset with the same size of the
    \textbf{CleanTestDataSET} that will be used to validate the performance of our model on
    the clean case after performing the tuning.
    
\subsubsection*{NoisyDataSET} Consisting of:\\
\indent $\bullet$ \textbf{NoisySmallTestDataSET}: again a dataset of
$N_{\text{test}}=7,500$ samples, but this time with noise. This is used to test
the robustness of the models trained by the \textbf{CleanTestDataSET};

$\bullet$ \textbf{NoisyTrainingDataSET}: a dataset with $N_{\text{train}}=150,000$
noisy samples that we use to train our models to be robust to noise;

$\bullet$ \textbf{NoisyTestDataSET}: a dataset with $N_{\text{test}}=50,000$ noisy
samples. This is the dataset that we use to test the performances of our models
trained with the \textbf{NoisyTrainingDataSET}.









\begin{figure}[hbpt]
  \pgfplotstableread{data/clean_and_noisy_datasets_svd.dat}{\svddata}%
  \centering
  \begin{subfigure}[b]{.495\textwidth}
    \tikzsetnextfilename{clean_dataset_svd}
    \begin{tikzpicture}
      \begin{axis}[
          width=.95\textwidth,
          height=.5\textwidth,
          axis x line=bottom,
          axis y line=left,
          xmax=200,
          ymin=1e-5, ymax=1e4,
          ymode=log,
          ylabel = {},
          xlabel = {}
        ]
        \addplot [blue, very thick, dashed] table[x expr = \coordindex + 1, y = clean-c-singular-values]{\svddata};
        \addplot [orange, very thick] table[x expr = \coordindex + 1, y = clean-u-singular-values]{\svddata};
      \end{axis}
    \end{tikzpicture}
  \end{subfigure}
  \begin{subfigure}[b]{.495\textwidth}
    \tikzsetnextfilename{noisy_dataset_svd}
    \begin{tikzpicture}
      \begin{axis}[
          width=.95\textwidth,
          height=.5\textwidth,
          axis x line=bottom,
          axis y line=left,
          xmax=400,
          ymin=1e-5, ymax=1e4,
          ymode=log,
          ylabel = {},
          xlabel = {}
        ]
        \addplot [blue, very thick, dashed] table[x expr = \coordindex + 1, y = noisy-c-singular-values]{\svddata};
        \addplot [orange, very thick] table[x expr = \coordindex + 1, y = noisy-u-singular-values]{\svddata};
      \end{axis}
    \end{tikzpicture}
  \end{subfigure}
  \caption{SVD analysis of \textbf{CleanTrainingDataSET} and \textbf{NoisyTrainingDataSET}. The
  figure shows the magnitude of the first 200 singular values for
  \textbf{CleanTrainingDataSET} (left) and the first 400 singular values for
  \textbf{NoisyTrainingDataSET}. The singular values for the doping profile matrices are
  shown in dashed blue lines and the singular values for the photovoltage
  signals are shown as solid orange lines.}
  \label{fig:svd-datasets}
\end{figure}
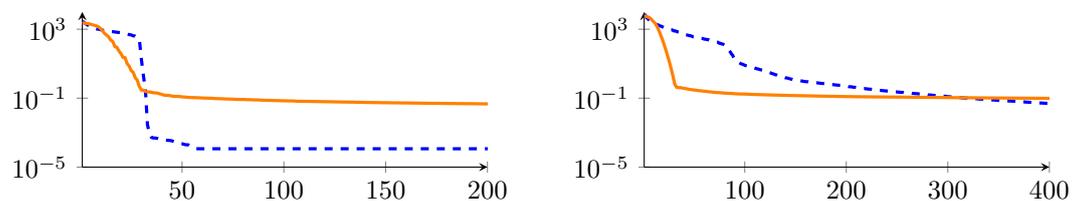

We provide a qualitative analysis of \textbf{CleanTrainingDataSET} and
\textbf{NoisyTrainingDataSET} in \Cref{fig:svd-datasets}, where we compute the
Singular Value Decomposition (SVD) of both datasets interpreted as $n\times
N_{\text{train}}$ matrices. The number of dominant singular values is a crude
indication of the actual dimension of the spaces $\DopingSpaceLPS$, and
$\SignalSpaceLPS$, and shows that, roughly, their dimension is close to $50$
when using the clean data generation described above, see
\Cref{fig:svd-datasets} left panel. The two dimensions mismatch significantly
when adding noise, see \Cref{fig:svd-datasets} right panel. While the
dimension of the clean photovoltage signals follows roughly that of the
clean doping profiles, the same cannot be said for the noisy cases. This
difference is a hint that the inverse operator $\Inverse$ is not well-posed, due
to a mismatch in the dimension of the input and output spaces.

\subsection{Least squares}
\label{ssec:res-linear}

In \cite{Farrell2021} the authors showed that, under certain restrictive conditions, the
operator $\Inverse$ can be approximated by a linear one. Indeed, in case of an
$n$ or $p$ doped semiconductor that varies only along the $x$ direction, we can show that
\begin{equation}
  u_P(\mathbf{x}, C)\sim
  -\mathbf{e}_x \cdot \nabla C(\mathbf{x}).
\label{eq:u_LPS-prop_N_D}
\end{equation}
If the support of the laser source is larger than the wavelength of oscillation 
of the doping profile, measuring the voltage difference may actually result in a signal that does not capture a single oscillation, but a floating average. 
We could represent this floating average of
doping fluctuations by a convolution of the doping gradient with some (unknown)
profile $f$, i.e.\
\begin{equation}
  u_P(x, C)\sim (f \ast( \mathbf{e}_x \cdot \nabla C))(x) = \dfrac{\mathrm{d}}{\mathrm{d}x}\left(f \ast C \right)(x).
\label{eq:convolutionWithGradient}
\end{equation}
Hence, it may seem plausible to relate the photovoltage signal to the doping
profile via \textit{linear} operations such as convolution and differentiation.
However, it is not clear whether this profile function $f$ is actually
independent from the doping fluctuations, the shape of the laser beam or the
laser spot position.

The least square analysis is useful to understand how the inverse problem
behaves, and what type of operation the inverse operator $\Inverse_h(\mathbf{u})
= \mathbf{A}_h \mathbf{u}$ performs. Let us express the training data for the
photovoltage signals and doping profiles generated according to
\Cref{sec:data_generation} with the two dense matrices $\mathbf{U}^{n\times
N_{\text{train}}}$ and $\mathbf{C}^{n\times N_{\text{train}}}$. Then  we wish to
solve the least square problem
\begin{equation}
  \label{eq:least-square-matrix}
  \mathbf{A}_h = \argmin_{\mathbf{B} \in \mathbb{R}^{n\times n}} \frac12 \left\| \mathbf{B} \mathbf{U}^{n\times N_{\text{train}}} - \mathbf{C}^{n\times
  N_{\text{train}}} \right\|^2_{\ell^2}.
\end{equation}

\begin{figure}
  \centering
  \includegraphics[width=.64\textwidth]{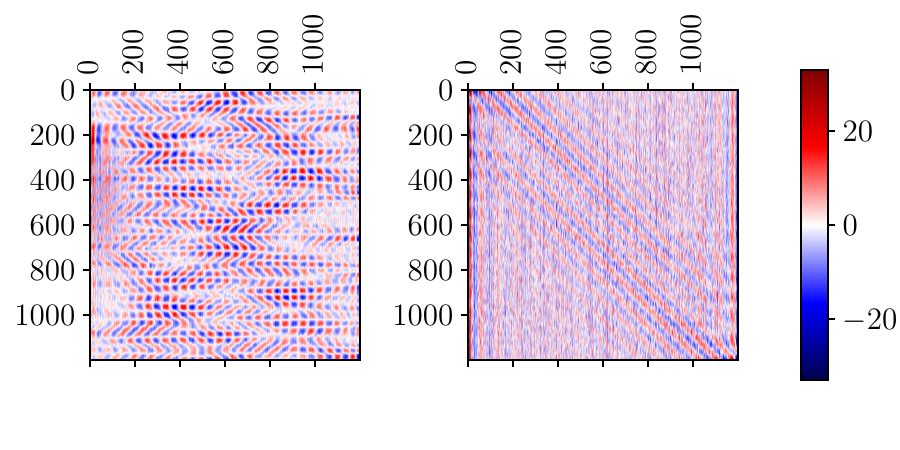}
  \hspace{-0.95em} 
  \tikzsetnextfilename{ls-matrices-svd-values}
  \begin{tikzpicture}
    \pgfplotstableread{data/ls_matrix_s_values.dat}{\lsSValues}%
    \begin{axis}[
        width=.33\textwidth,
        height=.33\textwidth,
        axis x line=bottom,
        axis y line=left,
        ymode = log,
        ymin =1e-1,
        ylabel = {},
        xlabel = {}
      ]
    \addplot [blue, dashed, very thick] table[x expr = \coordindex, y = clean]{\lsSValues};
    \addplot [orange, very thick] table[x expr = \coordindex, y = noisy]{\lsSValues};
    \end{axis}
  \end{tikzpicture}
  \caption{Magnitude of the entries of the least square matrices $\mathbf{A}_h$ obtained for 
  \textbf{CleanDataSET} (left) and \textbf{NoisyDataSET} (center), and singular values of the two matrices for the clean case (dashed blue line) and for the noisy case (orange line).}
  \label{fig:least-square-matrices}
\end{figure}
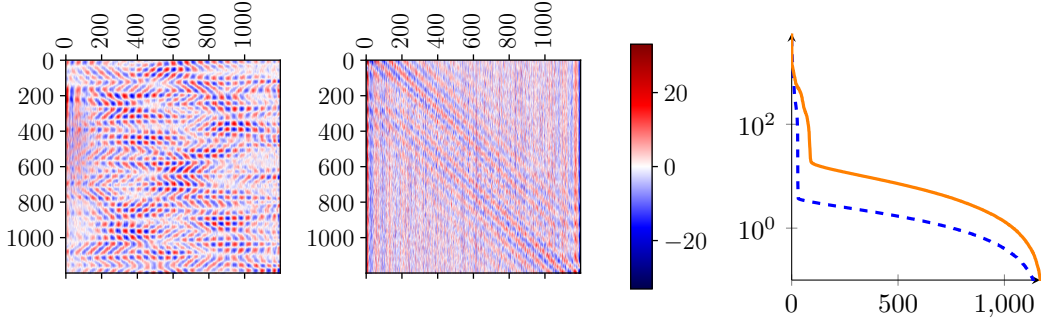

In \Cref{fig:least-square-matrices} we show $\mathbf{A}_h$ obtained with 
\textbf{CleanDataSET} (left) and with  \textbf{NoisyDataSET} (center). We observe a
highly non-local behavior of $\Inverse_h$. In the \textbf{CleanDataSET} case, this behaviour is more pronounced than in the \textbf{NoisyDataSET} case.

If we compare the singular values of the matrix
$\mathbf{A}_h$ in \Cref{fig:least-square-matrices} (right) with those of
the two datasets in \Cref{fig:svd-datasets}, we observe that the dominant singular values of the
least square matrix computed with the \textbf{CleanDataSET} are the first thirty
(similar to what happens in the singular values of the \textbf{CleanDataSET}
itself in \Cref{fig:svd-datasets}), while those that are most relevant for the
matrix constructed with the \textbf{NoisyDataSET} are the first one hundred.
While this information is only qualitative, it does show that a non-negligible part of the local inverse problem can be approximated relatively well by the linear
 operator $\Inverse_h(\mathbf{u}) = \mathbf{A}_h \mathbf{u}$, with an intrinsic dimension of around one
hundred, when including noise, and much smaller when considering a
clean dataset.

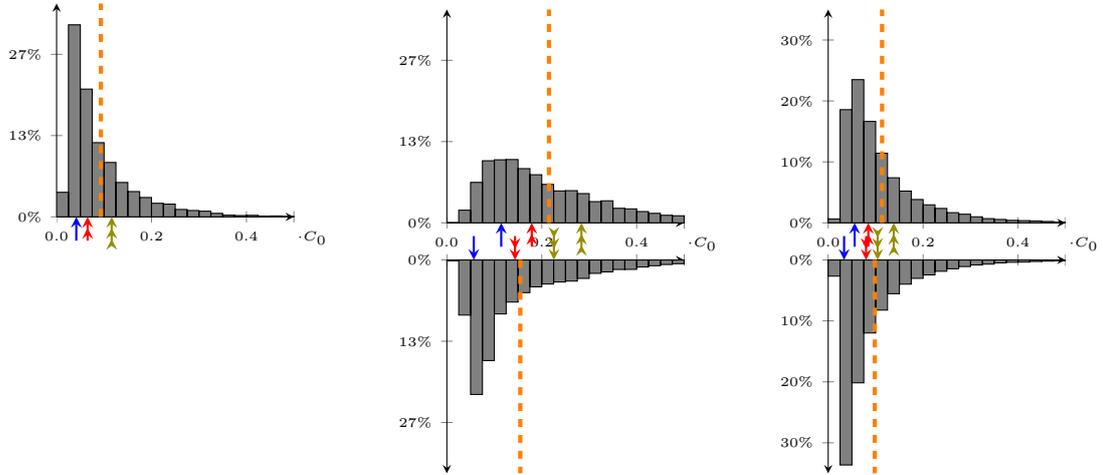
\begin{figure}
    \tikzsetnextfilename{ls-histogram}
    \begin{tikzpicture}
      \DrawHistogramPrecomputed[width=.32\textwidth, height=.3\textwidth, hbin=35,
      topmean=0.0932264165148, 
      npoints=7500]%
        [4.15e-2]
        [6.53e-2]
        [1.16e-1]
        {data/ls_linfty_absolute_errors_precomputed.csv}
    \end{tikzpicture} 
    \hfill
    \tikzsetnextfilename{ls-histogram-clean-noisy}
    \begin{tikzpicture}
      \DrawHistogramWithAndWithoutMeanPrecomputed[%
          width=.32\textwidth, height=.3\textwidth, 
          hbin=35,
          topmean=0.215282344165,
          bottommean=0.15471798347,
          npoints=7500%
        ]%
        [1.14565735598e-1]
        [5.67525588995e-2]
        [1.78950190432e-1]
        [1.43667801856e-1]
        [2.83614753204e-1]
        [2.25586040729e-1]
        {data/ls_linfty_absolute_errors_clean_noisy_precomputed.csv}
    \end{tikzpicture}
    \hfill
    \tikzsetnextfilename{ls-histogram-noisy-noisy}
    \begin{tikzpicture}
     \DrawHistogramWithAndWithoutMeanPrecomputed[%
         width=.32\textwidth, height=.3\textwidth, 
         hbin=35, 
         npoints=5e4,%
         topmean=1.13500970111e-1, 
         bottommean=9.80772295572e-2
       ]%
       [5.57298390205e-2]
       [3.35201578969e-2]
       [8.4565471851e-2]
       [7.97722105158e-2]
       [1.38122735046e-1]
       [1.04144876472e-1]
       {data/ls_linfty_absolute_errors_noisy_noisy_precomputed.csv}
   \end{tikzpicture}
   \caption{Statistical distribution of the errors on the predictions of our least
   square model, trained/tested on
   \textbf{CleanTrainingDataSET}/\textbf{CleanTestDataSET} (left),
   \textbf{CleanTrainingDataSET}/\textbf{NoisySmallTestDataSET} (center), and
   \textbf{NoisyTrainingDataSET}/\textbf{NoisyTestDataSET}. The bottom histograms
   are obtained by removing the average value of the doping. The dashed lines
   represent the average value of the error, while the arrows point to the 25, 50,
   and 75-percentile of the error on the top histogram, and show how these
   change when removing the average in the bottom histogram.}
   \label{fig:LS-histograms}
  \end{figure}

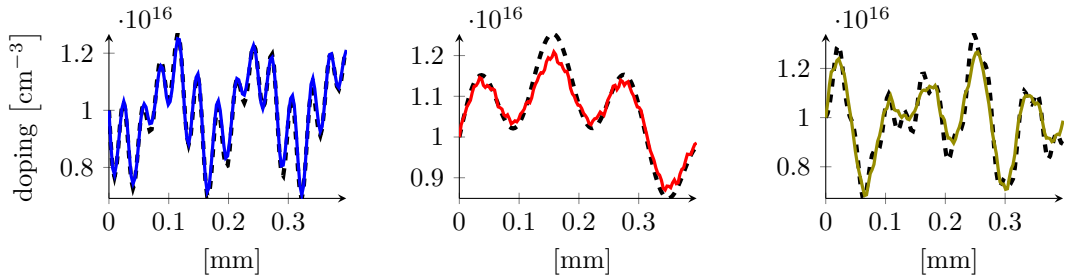
\begin{figure}
  \centering
  \begin{subfigure}[b]{.32\textwidth}
    \centering
    \tikzsetnextfilename{ls-example-p25-cleanclean}
    \begin{tikzpicture}
      \DrawSample[width=\textwidth, height=.8\textwidth, color=blue]%
           {data/ls_examples.csv}{p25doping}{p25prediction}
    \end{tikzpicture}
  \end{subfigure}
  \begin{subfigure}[b]{.32\textwidth}
    \centering
    \tikzsetnextfilename{ls-example-p50-cleanclean}
    \begin{tikzpicture}
      \DrawSampleNoLabel[width=\textwidth, height=.8\textwidth]%
           {data/ls_examples.csv}{p50doping}{p50prediction}
    \end{tikzpicture}
  \end{subfigure}
  \begin{subfigure}[b]{.32\textwidth}
    \centering
    \tikzsetnextfilename{ls-example-p75-cleanclean}
    \begin{tikzpicture}
      \DrawSample[width=\textwidth, height=.8\textwidth, color=olive]%
           {data/ls_examples.csv}{p75doping}{p75prediction}
    \end{tikzpicture}
  \end{subfigure}
 \caption{%
   Three examples of predictions obtained using the least squares model applied
   to the \textbf{CleanDataSET}. The dashed gray line is the expected result,
   the continuous colored lines are our predictions. The three plots represent
   samples whose error is close to the 25, 50 and 75-percentile (from left to
   right), and corresponds to the three arrows in the left histogram in
   \Cref{fig:LS-histograms} of the same color.}
 \label{fig:LS-examples}
\end{figure}

The statistical distribution of the test error in the \textbf{CleanDataSET} (i.e.,
the error defined in \eqref{eq:testing-error} for the \textbf{CleanTestDataSET})
is reported in \autoref{tab:errors} and depicted
in the left panel of \autoref{fig:LS-histograms}. 
It shows an error which is, on average (orange dashed line), around
$9\%$. Even though such a model presents a relatively high error, we observe
that it is still capable of capturing the overall qualitative behavior of the
doping profiles (see \Cref{fig:LS-examples}).

Next, we check how robust our two models are with respect to noise. As a first test, we use
the linear operator generated from the \textbf{CleanTrainingDataSET} to predict the doping of
the samples in the \textbf{NoisyTestDataSET}. Since this model has not been trained with
data affected by noise, we observe a significant deterioration in
the quality of the predictions, shown in the middle plot of  \Cref{fig:LS-histograms},
where the average error is now around $20\%$.

We believe that the main reason why the error increases so much is the fact that
the noise introduces a shift in the average of the doping function that can not
be physically estimated by our setup. In other words, 
\Cref{eq:doping_space} defines the admissible doping profiles with an 
average of (roughly) equal to $1\cdot 10^{16}=C_0$ on the
overall domain.
When
we introduce noise, this assumption does not hold anymore and we have no way
to predict whether the average of the doping is the one we expect in the entire
crystal sample. This is also coherent with the qualitative analysis in \Cref{eq:u_LPS-prop_N_D}, where we show that the value of the current is
mostly related to the local variation of the doping and not to its average value.

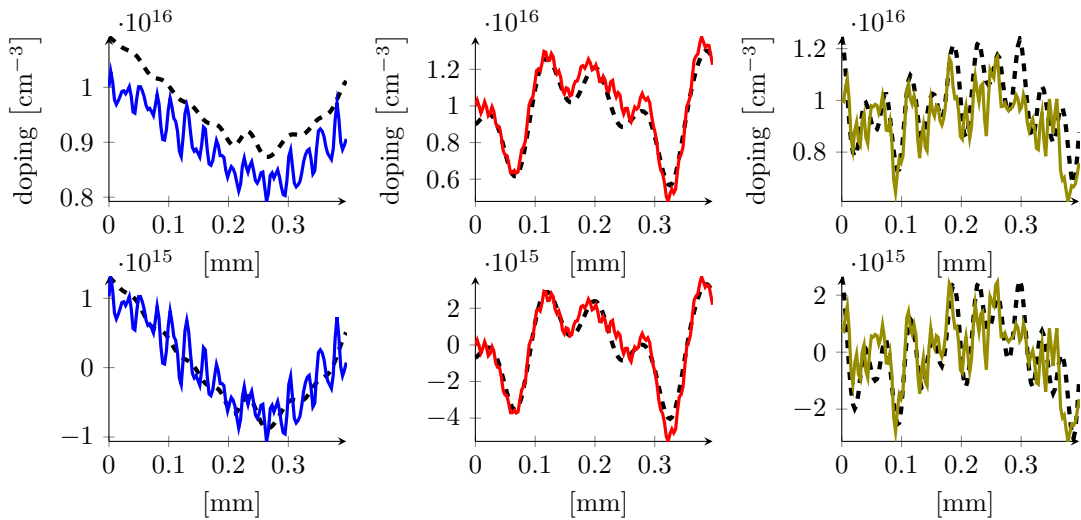
\begin{figure}
  \centering
 \begin{subfigure}[b]{.32\textwidth}
   \centering
   \tikzsetnextfilename{ls-example-p25-cleannoisy}
   \begin{tikzpicture}
     \DrawSampleWithAndWithoutMean*[width=\textwidth, height=1.6\textwidth, color=blue]%
          {data/ls_examples_clean_noisy.csv}{p25doping}{p25prediction}
   \end{tikzpicture}
 \end{subfigure}
 \begin{subfigure}[b]{.32\textwidth}
   \centering
   \tikzsetnextfilename{ls-example-p50-cleannoisy}
   \begin{tikzpicture}
     \DrawSampleWithAndWithoutMean*[width=\textwidth, height=1.6\textwidth]%
          {data/ls_examples_clean_noisy.csv}{p50doping}{p50prediction}
   \end{tikzpicture}
 \end{subfigure}
 \begin{subfigure}[b]{.32\textwidth}
   \centering
   \tikzsetnextfilename{ls-example-p75-cleannoisy}
   \begin{tikzpicture}
     \DrawSampleWithAndWithoutMean*[width=\textwidth, height=1.6\textwidth, color=olive]%
          {data/ls_examples_clean_noisy.csv}{p75doping}{p75prediction}
   \end{tikzpicture}
 \end{subfigure}
\caption{Three examples of predictions on the \textbf{NoisyTestDataSet}, using
the least square model trained on the \textbf{CleanTrainingDataSet}, taken from
the 25, 50, and 75 percentile of the error (from left to right), without removing the
average of the doping (top), and removing the average (bottom). The samples in the plots have an error that corresponds to the arrows in the histograms in the middle of \Cref{fig:LS-histograms} (including the color).}
\label{fig:LS-examples-clean-noisy}
\end{figure}

Removing the average from both the output of the model, and from the reference doping during testing leads to the results in the right plot of \Cref{fig:LS-histograms}
(bottom), where the error drops from $20\%$ to $15\%$. Indeed, we expect still a
larger error w.r.t. to the histogram presented in \Cref{fig:LS-histograms}, 
since the test dataset includes noisy 
data which are not included in the training dataset.
A few examples of predictions of the doping profile on $\Sigma_h$ with noisy input data are shown in
\Cref{fig:LS-examples-clean-noisy}.

An improvement of the error for the least square problem can be obtained by
performing the training on the \textbf{NoisyTrainingDataSET}. The histogram of
the error on the \textbf{NoisyTestDataSET} is shown in the right plot of
\Cref{fig:LS-histograms}, where the error is now around $10\%$ (see
\Cref{tab:errors} for the details). It becomes apparent that training with noisy
doping profiles  also significantly reduces the need to remove the average
doping value.


In summary, the least square model is able to
predict doping profiles with an average error of around $10\%$
for clean test data or noisy but properly average adjusted test data.
However, it is possible to further improve the results by introducing
nonlinearities in our models, for example using neural networks. We focus on multilayer perceptrons and residual neural networks in the
following sections.

\subsection{Multilayer perceptron}
\label{ssec:res-mlp}

While the least square approach is a good starting point, its efficiency is only
good if the inverse operator is close to a linear operator. For inverse problems
associated to the van Roosbroeck system, this is not necessarily the case, and
one may choose to introduce some nonlinearities in the approximation of
$\Inverse_h$.
Multilayer perceptrons \cite{Goodfellow2016} are among the most widely used
feedforward neural networks, and they are the first natural candidate for
general nonlinear function approximations.
We consider networks consisting of a down-sampler ($L_1$), a multilayer
perceptron with six layers ($L_2\hdots L_7$), and an up-sampler ($L_8$).
The main motivation to introduce the down-sampler and the up-sample layers is
to avoid introducing ``too many values'' inside our model. Indeed, if we used
directly the data from our datasets, then the first layer of our MLP should
have 1200 neurons (which is the value of $n$ defined in \cref{sec:setup}),
which is a large number that would increase the complexity of the model and
our probability of overfitting during the training. On the other hand, the
signal is just a spatial function defined on $\Sigma$ and we can interpolate
the signal on a coarser grid obtaining a still accurate description
of the signal while reducing the dimension of the space of the admissible
inputs. The down-sampler layer we introduced
uses cubic interpolation to describe the original signal on a coarse grid
of evenly spaced points. The up-sampler, instead, performs cubic
interpolation from the coarse grid to the original one we used for measuring
the doping. In this way, we can compare results from different models in the
same space, independently from the grid used by the MLP.

There is no a priori method to choose a suitable number of points for
the coarse grid nor there exists a standard algorithm to select the number
of neurons of each layer.
Our strategy is to choose some reasonable values for each free hyperparameter
of the model and then combine all the admissible choices to
generate a set of candidate models. Unfortunately, we do not dispose of enough
computational power to train all these models and, therefore, we need an algorithm
to explore this set and to find a good choice for the hyperparameters.

First of all, let us describe the space of the admissible configurations.
The number of points of the coarse grid (which coincides with the number
of neurons of the input layer of the MLP) is chosen from the set $\{100 + i 50: i =
0,\ldots, 8\}$. The number of neurons of the MLP output layer $L_7$, which in
principle can be different from the previous one, is chosen from the same set.

Let $\#(L_i)$ denote the number of neurons in $L_i$. For each hidden layer $L_i$
of the MLP with $i=3,\ldots,6$, we randomly choose $\#(L_i)$ in the set
$\{\#(L_{i-1}),\#(L_{i-1})\pm 50, \#(L_{i-1})\pm 100, \#(L_{i-1})-200\}$ with
the constraint that $\#(L_i)>0$. In total, there are $71,118$ admissible
configurations. There are fewer than $6^4\cdot 9^2$ configurations because some of them would lead to a negative amount of neurons. 
The multilayer perceptrons are implemented in PyTorch \cite{NEURIPS2019_9015}.


For a given configuration of the neural network, we apply the stochastic
gradient descent (SGD) algorithm (without momentum and weight decay) to find the
MLP that minimizes the mean square error loss defined in \Cref{eq:msel}. For
each network configuration, we fix the learning rate to a value that is chosen
randomly in the interval $[10^{-3}, 1]$ and the batch size to $64$. The number of
training epochs is set to $1,000$. 

We randomly select $10,000$ configurations where we vary both the structure of
the MLP (choosing one of the 71,118 possible MLPs) as well as the learning
parameters of SGD. To reduce the computational cost associated with the training
of all the resulting neural networks, we apply the Asynchronous Successive
Halving algorithm (ASHA)~\cite{ashaalgorithm}, a multi-armed bandit algorithm
that has been optimized to run on a massive amount of parallel machines. 
The ASHA algorithm discards the worst $50\%$ performers of the neural network
population at each check point, and proceeds with the training only for the most
promising neural networks. Our implementation is based on the Ray
library~\cite{ray}, which also takes care of distributing the jobs and
coordinating the execution among different nodes during the computation.



\begin{figure}
  \tikzsetnextfilename{mlp-histogram}
  \begin{tikzpicture}
   \DrawHistogramPrecomputed[%
       width=.32\textwidth, 
       height=.3\textwidth, 
       hbin=65,
       npoints=7500,
       topmean=0.0466581875152%
     ]%
     [1.03e-2]
     [1.83e-2]
     [3.98e-2]
     {data/mlp_linfty_absolute_errors_precomputed.csv}
  \end{tikzpicture}
  \hfill
  \tikzsetnextfilename{mlp-histogram-clean-noisy}
  \begin{tikzpicture}
    \DrawHistogramWithAndWithoutMeanPrecomputed[%
        width=.32\textwidth, height=.3\textwidth,
        hbin=65, npoints=7500, topmean=0.211890530963,
        bottommean=0.137045126794%
      ]%
      [9.10514550976e-2]
      [8.48619686971e-2]
      [1.73469650447e-1]
      [5.7381736591e-2]
      [2.92739070476e-1]
      [1.99446428532e-1]
      {data/mlp_linfty_absolute_errors_clean_noisy_precomputed.csv}
  \end{tikzpicture}
  \hfill
  \tikzsetnextfilename{mlp-histogram-noisy-noisy}
  \begin{tikzpicture}
    \DrawHistogramWithAndWithoutMeanPrecomputed[%
        width=.32\textwidth, height=.3\textwidth, hbin=65, npoints=5e4,%
        topmean=7.49506840602e-2, bottommean=5.39628229523e-2%
      ]%
      [2.13785874311e-2]
      [1.76963068753e-2]
      [3.59810025679e-2]
      [2.11522834264e-2]
      [6.20446002489e-2]
      [4.75789006984e-2]
      {data/mlp_linfty_absolute_errors_noisy_noisy_precomputed.csv}
  \end{tikzpicture}
  \caption{Statistical distribution of the $\ell^\infty$ errors on the predictions of our best
 multilayer perceptron, trained/tested on
 \textbf{CleanTrainingDataSET}/\textbf{CleanValidationDataSET} (left picture),
 \textbf{CleanTrainingDataSET}/\textbf{NoisySmallTestDataSET} (center), and
 \textbf{NoisyTrainingDataSET}/\textbf{NoisyTestDataSET} (right). The bottom histograms
 are obtained by removing the average value of the doping. The dashed lines
 represent the average value of the error, while the arrows point to the 25, 50,
 and 75-percentile of the error on the top histogram, and show how these
 transform when removing the average in the bottom histogram.}
  \label{fig:MLP-histograms}
 \end{figure}
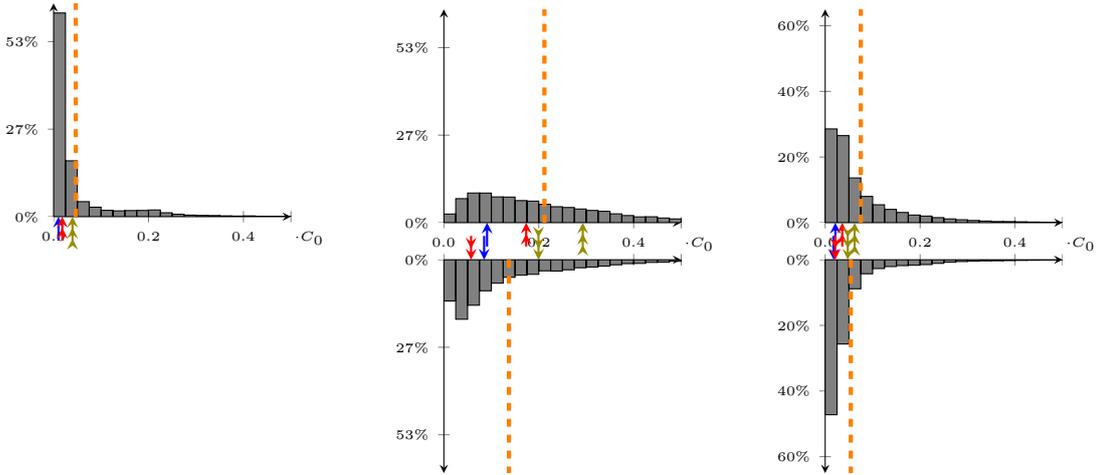
 
 The multilayer perceptron with the smallest $\ell^\infty$ testing error defined
in \Cref{eq:testing-error} (trained with \textbf{CleanDataSET} and tested with
\textbf{CleanTestDataSET}), has six layers with $250, 250, 150, 100, 1000, 350$
nodes, a batch size of $64$, and a learning rate of $0.06$. The statistical
distribution of the error is reported in the middle part of \Cref{tab:errors}, and depicted in the
left plot of \autoref{fig:MLP-histograms}. The $\ell^\infty$ error is, on
average, around $4.67\%$. Even with such a simple neural network, we reduce the
average error with respect to the least square model by of a factor two. Also,
we obtain a much better statistical distribution, clearly shown in
\Cref{fig:MLP-histograms}.


Furthermore, the resulting MLP is robust to noise, provided that we properly
filter the average of the doping profile. A comparison between the statistical
distribution of the errors (both with and without average) as well as the
corresponding statistical distributions are shown in \Cref{tab:errors} and in
\Cref{fig:MLP-histograms} (center). Testing with the \textbf{NoisyTestDataSET}
while the MLP was trained with \textbf{CleanDataSET} shows (as expected) a large
increase in the average error.

\begin{figure}
   \centering
  \begin{subfigure}[b]{.32\textwidth}
    \centering
    \tikzsetnextfilename{mlp-example-p25-noisynoisy}
    \begin{tikzpicture}
      \DrawSampleWithoutMean[width=\textwidth, height=0.8\textwidth, color=blue]%
           {data/mlp_examples_noisy_noisy.csv}{p25doping}{p25prediction}
    \end{tikzpicture}
  \end{subfigure}
  \begin{subfigure}[b]{.32\textwidth}
    \centering
    \tikzsetnextfilename{mlp-example-p50-noisynoisy}
    \begin{tikzpicture}
      \DrawSampleWithoutMean[width=\textwidth, height=0.8\textwidth]%
           {data/mlp_examples_noisy_noisy.csv}{p50doping}{p50prediction}
    \end{tikzpicture}
  \end{subfigure}
  \begin{subfigure}[b]{.32\textwidth}
    \centering
    \tikzsetnextfilename{mlp-example-p75-noisynoisy}
    \begin{tikzpicture}
      \DrawSampleWithoutMean[width=\textwidth, height=0.8\textwidth, color=olive]%
           {data/mlp_examples_noisy_noisy.csv}{p75doping}{p75prediction}
    \end{tikzpicture}
  \end{subfigure}
 \caption{Three examples of predictions on the \textbf{NoisyTestDataSet}, with
 the best MLP model trained on the \textbf{NoisyTrainingDataSet}, taken from the
 25, 50, and 75 percentile of the error, removing the average. The errors of the plots correspond to the arrows in the lower histogram in the right of \Cref{fig:MLP-histograms}.}
 \label{fig:MLP-examples-noisy-noisy}
\end{figure}
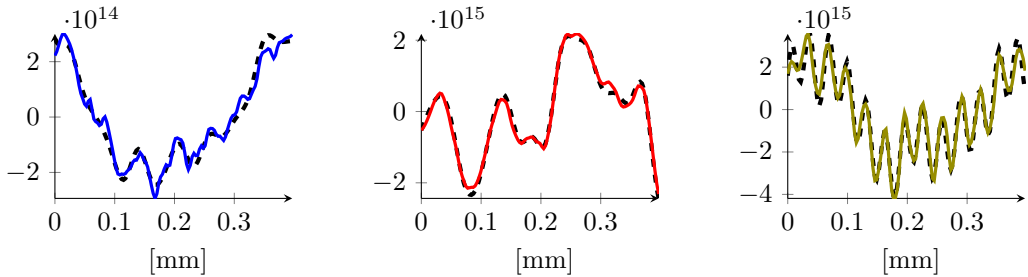

More robust results with respect to noise can be obtained by training again the
best MLP on the full \textbf{NoisyTrainingDataSET}. We do not perform an
additional hyperparameter tuning, but simply repeat the training stage on the
same MLP.
In this case, the larger dataset and the presence of noise induce an error when
tested with the \textbf{NoisyTestDataSET} of around $5.92\%$. This error can,
again, be improved by removing the average as we did before, leading to a final
error which is very close to the clean case ($4.96\%$ vs $4.67\%$). We show three doping reconstructions for this case in the 25,
50, and 75 percentile in \Cref{fig:MLP-examples-noisy-noisy}.

In summary, introducing a nonlinear MLP model has reduced the error
roughly by a factor of two compared to the linear least square model. The average MLP is just below $5\%$
for clean test data or noisy but properly average adjusted test data.

\subsection{Residual neural network (ResNet)}
\label{ssec:res-ResNet}

In 2015, Kaiming et al.\cite{KaimingZhang2015} introduced ResNets,
feedforward convolutional neural network models, which, since then, have been widely
used to solve different kinds of problems.
One of the biggest advantages of ResNets (or, generally, of convolutional neural networks) over the multilayer perceptrons is their ability to learn 2D images. Interpreting the doping profiles in 2D probing region as an image, this is precisely what we need to achieve, too. 
Moreover, the proportionality shown in  \Cref{eq:convolutionWithGradient}, valid under additional assumptions, suggests that 
the photovoltage signal and the doping profile may be related
by convolutions. Even assuming now a nonlinear relationship described via ResNets, we may still wish to exploit a convolutive structure to encompass the relationship in  \Cref{eq:convolutionWithGradient}.
Our reference implementation for residual neural networks (ResNet) is described
in \cite{KaimingZhang2015}, with some important differences regarding the
structure of the network.

First, in our LPS setup from \Cref{sec:setup}, we consider one-dimensional data vectors $\mathbf{C}$ and $\mathbf{u}$ instead of two-dimensional images. The dimension is particularly
 relevant to choose the size of the network: In ResNets, the downsampling
 blocks reduce by a factor of two the size of the input, by first halving the size
 of the input in each direction (reducing the spatial indices by a factor four) and then
 doubling the number of channels. In our setup, this dimension reduction
 associated to the downsampling block does not happen because in the first step
 we only have one dimension to divide by two, and we end up with the same number
 of neurons after doubling the numbers of channels.
Moreover, while the ImageNet database used to train the network described in
\cite{KaimingZhang2015} contains about one million images, our datasets are significantly
smaller. We expect to face some overfitting and therefore we aim for a model
with fewer parameters than the one described in
\cite{KaimingZhang2015}.
Finally, we solve a regression problem instead of a classification one. The
major consequences for our model are that we cannot use drop-out layers to
reduce overfitting, and it is unclear if there is any benefit in using batch
normalization layers. We keep the batch normalization layers because we
expect the statistical distribution of our batches to be similar to each other,
and the statistical distribution of the error for a simpler network (the MLP)
shows an improvement in the performance of the model with noise when removing
the average, suggesting that batch normalization layers are at least not
harmful.

We build several ResNet models using PyTorch (\cite{NEURIPS2019_9015}), and 
develop a strategy to find the best ones by adapting the model described in
\cite{KaimingZhang2015}, and using the insight obtained from the MLP model
presented in \Cref{ssec:res-mlp}.
Our best multilayer perceptron had 160,950 parameters that had to be optimized. We
try to keep the number of parameters of our ResNet around the same order of
magnitude. Of course, this means that we need to drastically simplify the model
developed in \cite{KaimingZhang2015}, to find a model whose number of parameters
is compatible with the size of our \textbf{CleanTrainingDataSET}.

Since training a ResNet is significantly more expensive compared to training a multilayer perceptron, we use a two-step strategy to reduce the computational cost of the hyperparameter tuning.
The possible configurations of the ResNet that we allow in our models are detailed in 
\Cref{sec:resnet-details}, and correspond to a total of $48$ configurations for
the gate, $9$ configurations for the encoder, and $18$ configurations for the
decoder, for a total of $\totalResNetConfigs$ possible configurations.

\begin{figure}[hbt]
  \centering
  \tikzsetnextfilename{resnet-lr-vs-error}
  \begin{tikzpicture}
    \DrawResnetScatterPlot[width=\textwidth, height=.3\textwidth]%
        {data/resnet_training_results.dat}
  \end{tikzpicture}
  \centering
  \tikzsetnextfilename{resnet-bsize-vs-error}
  \begin{tikzpicture}
    \DrawResnetTrainingVsBatches[width=.49\textwidth, height=.25\textwidth]%
        {data/resnet_training_vs_batch_size.dat}
  \end{tikzpicture}
  \hfill
  \tikzsetnextfilename{resnet-wdecay-vs-error}
  \begin{tikzpicture}
    \DrawResnetTrainingVsWeightDecay[width=.49\textwidth, height=.25\textwidth]%
        {data/resnet_training_vs_wdecay.dat}
  \end{tikzpicture}
  \caption{%
    $\ell^\infty$ errors of our model with respect to the parameters of the optimizer. In the top plot, the colors correspond to the networks RN1-RN4,
    as seen in \Cref{tab:resnet-models}. The  other two plots (bottom row) are batch size and weight decay density plots computed via \texttt{KernelDensity} of \texttt{sklearn.neighbors}
    with a Gaussian kernel and a bandwith of 0.02. Darker line colors correspond to bigger parameter values (batch size or weight decay).
  }
\end{figure}
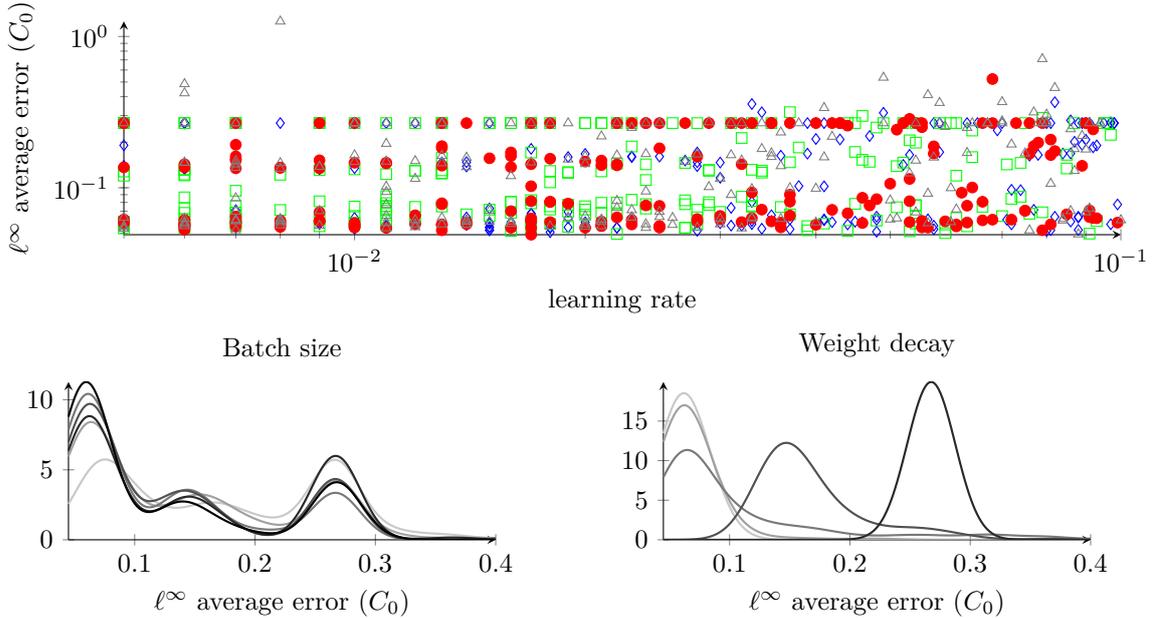

In the first step of our hyperparameter tuning, we restrict the parameters for the SGD optimizer to a fixed
float number chosen in the interval $[5\times 10^{-3}, 1\times 10^{-1}]$ for the  learning rate, and choose 
a batch size in the set $\{64, 128, 256, 384, 512\}$. We fix also
the gradient clipping to 1 and the weight decay to 0.
During this step, we randomly sampled 300 ResNet configurations among the \totalResNetConfigs generated according to \Cref{sec:resnet-details}, and for each one of them, we perform six trainings using different set of parameters for the SGD optimizer. Again, we rely on the ASHA algorithm as
we did in \Cref{ssec:res-mlp}. We compare the performances of the
models after 5000 epochs of the SGD algorithm, and retain only four of the candidate ResNet models selected by the ASHA algorithm, which are described in \Cref{tab:resnet-models}.

\begin{table}[h!tb]
  \centering
  \begin{tabular}{lcccc}
    & \textcolor{blue}{RN1} & \textcolor{green}{RN2} & %
      \textcolor{red}{\textbf{RN3}} & \textcolor{gray}{RN4} \\\hline
    N. of channels         & 24 & 16 & \textbf{24} & 16 \\
    Gate conv. kernel size &  9 &  5 &  \textbf{7} &  3 \\
    Gate conv. stride      &  4 &  4 &  \textbf{4} &  4 \\
    N. of encoder blocks   &  3 &  3 &  \textbf{3} &  3 \\
    Block type:            & \BlkFCB & \BlkBasic & \textbf{\BlkFCB} & \BlkFCB \\
    Downsampling:          & \TRUE & \TRUE & \textbf{\TRUE} & \TRUE \\
    Size of decoder layers & (100, 200) & (150, 150) & \textbf{(200, 200)} & (250, 200) \\
    N. of parameters       & 102,188 & 324,048 & \textbf{141,440} & 138,994
  \end{tabular}
  \caption{The four ResNet models (RN1--RN4) identified in the first stage of the hyperparameter tuning by the ASHA algorithm. RN3 is the one selected in the second stage.}
  \label{tab:resnet-models}
\end{table}

In the second stage of the hyperparameter tuning, we keep the four structures of the ResNet models fixed, and  enlarge the search space in the optimization parameters by defining a new set of admissible parameters: we choose the learning rate in the interval $[10^{-3}, 10^{-1}]$, the batch size
in the set $\{32, 64, 128, 256, 384, 512\}$, the gradient clipping in $\{10^{-2},
10^{-1}, 1\}$, and we introduce some weight decay, with decay parameters
chosen in the set $\{10^{-4}, 10^{-3}, 10^{-2}, 10^{-1}, 1\}$.

For each one of the four models, we perform 200 different trainings applying different parameters of the SGD optimizer, and finally select a winner, corresponding to RN3 in \Cref{tab:resnet-models}.

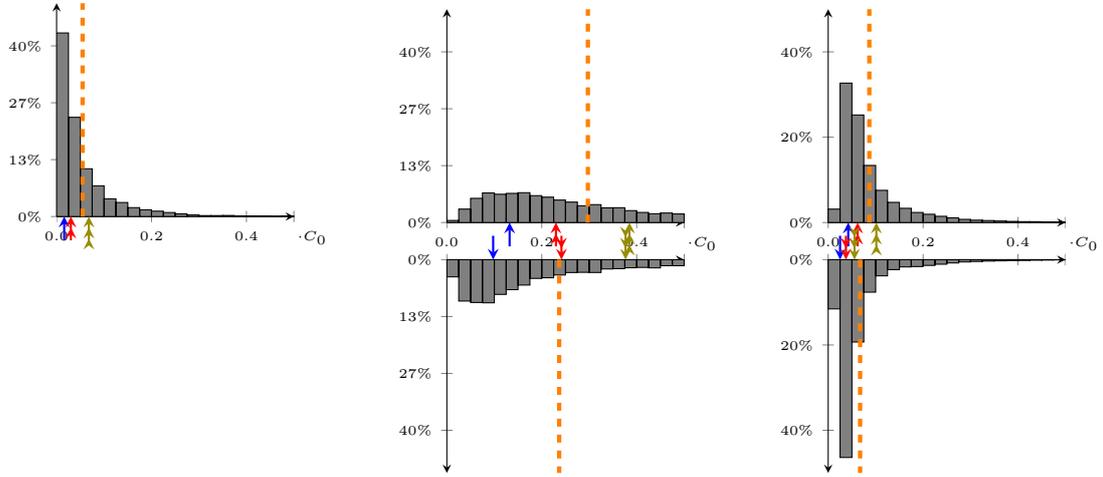
\begin{figure}
  \tikzsetnextfilename{resnet-histogram}
  \begin{tikzpicture}
   \DrawHistogramPrecomputed[%
       width=.32\textwidth, 
       height=.3\textwidth, 
       hbin=50,
       npoints=7500,
       topmean=0.0547251014013%
     ]%
     [1.60e-2]
     [3.00e-2]
     [6.80e-2]
     {data/resnet_linfty_absolute_errors_precomputed.csv}
  \end{tikzpicture}
  \hfill
  \tikzsetnextfilename{resnet-histogram-clean-noisy}
  \begin{tikzpicture}
    \DrawHistogramWithAndWithoutMeanPrecomputed[%
        width=.32\textwidth, height=.3\textwidth,
        hbin=50, npoints=7500, topmean=0.297174155529,
        bottommean=0.236715867679%
      ]%
      [1.32199268633e-1]
      [9.75741474942e-2]
      [2.29871549985e-1]
      [2.41478975152e-1]
      [3.84682294595e-1]
      [3.76539995839e-1]
      {data/resnet_linfty_absolute_errors_clean_noisy_precomputed.csv}
  \end{tikzpicture}
  \hfill
  \tikzsetnextfilename{resnet-histogram-noisy-noisy}
  \begin{tikzpicture}
    \DrawHistogramWithAndWithoutMeanPrecomputed[%
        width=.32\textwidth, height=.3\textwidth, hbin=50, npoints=5e4,%
        topmean=8.66750310472e-2, bottommean=6.72270677228e-2%
      ]%
      [4.22052024329e-2]
      [2.50142405781e-2]
      [6.20380418819e-2]
      [3.74011412158e-2]
      [1.0146075158e-1]
      [5.54264367475e-2]
      {data/resnet_linfty_absolute_errors_noisy_noisy_precomputed.csv}
  \end{tikzpicture}
  \caption{Statistical distribution of the errors on the predictions of our best
  ResNet, trained/tested on
  \textbf{CleanTrainingDataSET}/\textbf{CleanTestDataSET} (left),
  \textbf{CleanTrainingDataSET}/\textbf{NoisySmallTestDataSET} (center), and
  \textbf{NoisyTrainingDataSET}/\textbf{NoisyTestDataSET} (right). The bottom histograms
  are obtained by removing the average value of the doping. The dashed lines
  represent the average value of the error, while the arrows point to the 25,
  50, and 75-percentile of the error on the top histogram, and show how these
  transform when removing the average in the bottom histogram.}
  \label{fig:ResNet-histograms}
 \end{figure}

In \Cref{fig:ResNet-histograms} we present the statistical
distribution of the error for ResNet RN3 selected with the strategy outlined
above. We observe that the average error we obtain using the
\textbf{CleanDataSET} for both training and testing is around $5.47\%$ (see
\Cref{tab:errors} for the details). This number is slightly worse compared to the
MLP case introduced in \Cref{ssec:res-mlp}. In particular, we observe that the
ResNet RN3 trained with the \textbf{CleanDataSET} is more sensitive to
noise (see \Cref{fig:ResNet-histograms}, center figure) when compared to the
MLP ($29.7\%$ vs $21.2\%$). This difference decreases when removing the average, but for the ResNet remains surprisingly worse than the least square model ($16\%$ vs $15.5\%$).

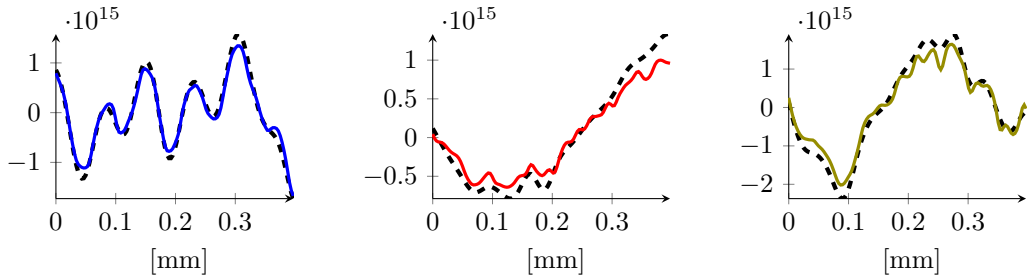
\begin{figure}
  \centering
 \begin{subfigure}[b]{.32\textwidth}
   \centering
   \tikzsetnextfilename{resnet-example-p25-noisynoisy}
   \begin{tikzpicture}
     \DrawSampleWithoutMean[width=\textwidth, height=0.8\textwidth, color=blue]%
          {data/resnet_examples_noisy_noisy.csv}{p25doping}{p25prediction}
   \end{tikzpicture}
 \end{subfigure}
 \begin{subfigure}[b]{.32\textwidth}
   \centering
   \tikzsetnextfilename{resnet-example-p50-noisynoisy}
   \begin{tikzpicture}
     \DrawSampleWithoutMean[width=\textwidth, height=0.8\textwidth]%
          {data/resnet_examples_noisy_noisy.csv}{p50doping}{p50prediction}
   \end{tikzpicture}
 \end{subfigure}
 \begin{subfigure}[b]{.32\textwidth}
   \centering
   \tikzsetnextfilename{resnet-example-p75-noisynoisy}
   \begin{tikzpicture}
     \DrawSampleWithoutMean[width=\textwidth, height=0.8\textwidth, color=olive]%
          {data/resnet_examples_noisy_noisy.csv}{p75doping}{p75prediction}
   \end{tikzpicture}
 \end{subfigure}
\caption{Three examples of predictions on the \textbf{NoisyValidationDataSet}, with
the best ResNet RN3 model trained on the \textbf{NoisyTrainingDataSet}, taken
from the 25, 50, and 75 percentile of the error, removing the average of the
doping. The colors of the plots indicate the arrows in the right bottom
histogram in \Cref{fig:ResNet-histograms}.}
\label{fig:ResNet-examples-noisy-noisy}
\end{figure}

When we keep the structure of RN3, and train the network on the \textbf{NoisyTrainingDataSET}, we obtain the results shown in \Cref{fig:ResNet-histograms} on the right. An example of reconstruction obtained with RN3 with a sample from the \textbf{NoisyTestDataSET} is shown in \Cref{fig:ResNet-examples-noisy-noisy}. In this case the average error on the \textbf{NoisyTestDataSET} is around $6.72\%$, which is slightly worse than the result obtained with the MLP network.

\section{Summary and outlook}
\label{sec:summary}

The aim of this paper was to properly model and numerically solve general
ill-posed inverse photovoltage technologies where measured photovoltage signals
are used to reconstruct local doping fluctuations in a semiconductor crystal.
The underlying charge transport model is based on the van Roosbroeck system as
well as  Ohm's law. 
We presented three different data-driven approaches to solve a physically relevant 
local inverse problem, namely via 
least squares, multilayer perceptrons, and residual neural networks.
The methods were trained on synthetic datasets (pairs of discrete doping
profiles and corresponding photovoltage signals at different illumination
positions) which are generated by efficient physics-preserving finite volume
solutions of the forward problem. 
While the linear least square
method yields an average absolute $\ell^\infty$ error of $9.3\%$, the nonlinear
networks roughly halve this error to $4.7\%$ and $5.5\%$, respectively, after
optimizing relevant hyperparameters. Our method turned out to be robust with
respect to noise, provided that training is repeated with larger, noisy, datasets. In this case, the error is around $9.8\%$, $5\%$, and $6.7\%$ respectively. Removing the average doping value from the data was more important to reduce the testing error for clean datasets than for noisy ones. The datasets, python codes, and resulting trained networks are available in the repository~\cite{ZenodoRepo2022}.

\begin{table}[ht]
  \begin{center}
    \begin{tabular}{c|l|cccc}
      \multicolumn{2}{c|}{Absolute error ($\times C_0$)}
       & Average & 25-th percentile & Median & 75-th percentile \\ \hline
      \multirow{5}{*}{LS}
       & \textbf{C.C.}          & $9.32\%$ & $4.15\%$ & $6.53\%$ & $1.16\%$ \\
       & \textbf{C.N.} w/ mean  & $21.5\%$ & $11.5\%$ & $17.9\%$ & $28.4\%$ \\
       & \textbf{C.N.} w/o mean & $15.5\%$ & $6.75\%$ & $10.5\%$ & $21.0\%$ \\
       & \textbf{N.N.} w/ mean  & $11.4\%$ & $5.72\%$ & $8.46\%$ & $13.8\%$ \\
       & \textbf{N.N.} w/o mean & $9.81\%$ & $4.20\%$ & $6.47\%$ & $11.9\%$ \\ \hline
      \multirow{5}{*}{MLP}
       & \textbf{C.C.}          & $4.67\%$ & $1.03\%$ & $1.83\%$ & $3.98\%$ \\
       & \textbf{C.N.} w/ mean  & $21.2\%$ & $91.1\%$ & $17.3\%$ & $29.3\%$ \\
       & \textbf{C.N.} w/o mean & $13.7\%$ & $4.15\%$ & $8.89\%$ & $19.7\%$ \\
       & \textbf{N.N.} w/ mean  & $5.92\%$ & $2.14\%$ & $3.60\%$ & $6.20\%$ \\
       & \textbf{N.N.} w/o mean & $4.96\%$ & $1.44\%$ & $2.49\%$ & $4.67\%$ \\ \hline
      \multirow{5}{*}{ResNet}
       & \textbf{C.C.}          & $5.47\%$ & $1.60\%$ & $3.00\%$ & $6.80\%$ \\
       & \textbf{C.N.} w/ mean  & $29.7\%$ & $13.2\%$ & $22.3\%$ & $38.5\%$ \\
       & \textbf{C.N.} w/o mean & $16.0\%$ & $5.27\%$ & $10.4\%$ & $21.2\%$ \\
       & \textbf{N.N.} w/ mean  & $8.67\%$ & $4.22\%$ & $6.20\%$ & $10.1\%$ \\
       & \textbf{N.N.} w/o mean & $6.72\%$ & $3.15\%$ & $4.44\%$ & $7.06\%$ \\ \hline
    \end{tabular}
  \end{center}
  \caption{Absolute errors w.r.t $C_0={1.0\times16}~{cm^{-3}}$ corresponding to the colored arrows as well as dashed lines showing in Figure~\ref{fig:LS-histograms}, \ref{fig:MLP-histograms}, and \ref{fig:ResNet-histograms} (from top to bottom). Here the second column categorized the training and testing datasets mentioned in the captions of the histogram figures. For instance, the row for ``\textbf{C.N.} w/ mean'' in the block of MLP shows the errors in the top-mid histogram in Figure~\ref{fig:MLP-histograms} which is trained/tested on \textbf{CleanTraining DataSET/NoisySmallTestDataSET}.}
  \label{tab:errors}
\end{table}

\begin{figure}
  \centering
  \tikzsetnextfilename{mlp-example-outlier}
  \begin{tikzpicture}
    \DrawSample[width=.49\textwidth, height=.3\textwidth, color=blue, colsep=comma]%
               {data/mlp_old_examples.csv}{d3}{p3}
  \end{tikzpicture}
   \tikzsetnextfilename{resnet-example-outlier}
   \begin{tikzpicture}
    \DrawSample[width=.49\textwidth, height=.3\textwidth, color=red]%
               {data/resnet_outlier.dat}{doping}{prediction}
  \end{tikzpicture}
  \caption{An example of an outliers using the MLP model (left) and the ResNet model RN3 (right). The dashed line represents the expected doping profile sample, and the continuous line is the output of the neural networks.}
  \label{fig:outliers}
\end{figure}
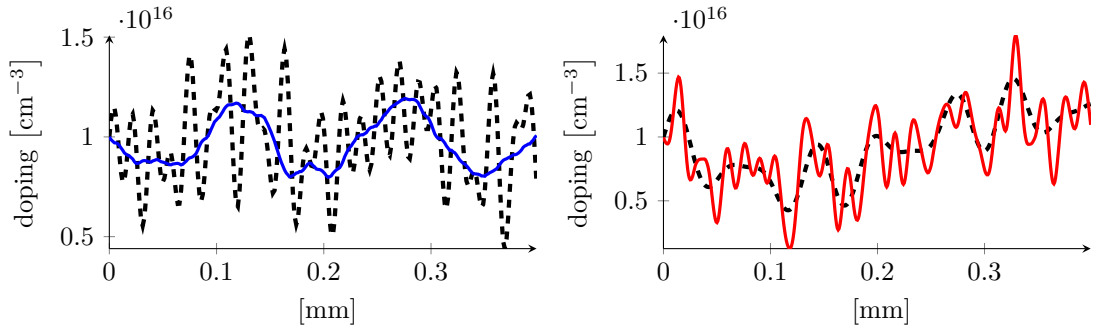

Future research may go into different directions: the numerical simulation
served as a proof of principle and was limited to a 2D version of the
photovoltage problem. For the 3D problem, efficient data generation strategies
need to be designed. 
Furthermore, it is clear that different doping profiles may
correspond to similar signals. On the one hand, this is intrinsically dependent on the technology that we apply to recover the signal: we hit the crystal with a laser that has a finite laser spot radius (in our case, around $20\mu m$), and we cannot expect to resolve oscillations with smaller wave lengths. On the other hand, our models do not have a way to distinguish between two doping profiles that deliver the same signal, but are exposed, during training, to many such cases. In \Cref{fig:outliers}, for example, we show two examples of doping profiles that have errors on the far right regions of the histograms in \Cref{fig:MLP-histograms} and \Cref{fig:ResNet-histograms}. The figure shows how in some cases, the inverse problem is oblivious to higher oscillations, and both neural network types will have cases in which they return an answer which is either too smooth w.r.t. the expected result, (left panel of \Cref{fig:outliers}) or too oscillatory (right panel of \Cref{fig:outliers}).

A possible focus of future research could be
to study how to resolve this ambiguity, which is intrinsic to the ill-posedness
of the discrete local photovoltage inverse problem.
Finally, we aim to extend our data-driven approach to opto-accoustic imaging, replacing the charge transport model with appropriate thermal expansion models.

\appendix
\section{Generation of noise}
\label{sec:noise}
In order to introduce noise in the doping profile \eqref{eq:doping_space}, we
perturb both the doping amplitudes as well as the wavelength. To perturb the
doping amplitudes, we start from some elements of $\DopingSpaceLPS$. We
introduce imperfections of the doping profile that may be interpreted physically
with a slight perturbation during the growth process (for example, due to a
fluctuation of the temperature), we choose a \textit{random} function $f_n(x)$
and define a perturbed doping of the form $\tilde{C}(x) \Def C(x) + f_n(x)$. We require that $f_n(x) \ll C(x)$ for any $x$. To generate $f_n(x)$, we pick 129
equally-spaced points $x_i$ across the silicon sample and randomly sample 129
values from a normal distribution (with 0 mean and standard deviation equal to 1),
obtaining a family of values $s_i$. Denoting the maximal variation of $C$ with
\[ 
    \Delta_C \Def \max_{x\in  [-\nicefrac{\ell}{2}, \nicefrac{\ell}{2}]}C(x) -
                  \min_{x\in  [-\nicefrac{\ell}{2}, \nicefrac{\ell}{2}]}C(x)
\quad 
\text{we define} \quad 
    f_n(x_i) \Def k \,s_i \, \Delta_C, 
\]
for $0<k \ll 1$. For any other point $x\neq x_i$ across the sample, we compute  $f_n(x)$ by cubic
spline interpolation. Next, we perturb the wavelengths. In \eqref{eq:doping_space},
we assumed that the doping has fluctuations with constant wavelength across the
entire domain. We weaken this assumption by introducing a non periodic
perturbation in the argument of the sinusoidal functions: we consider a
differentiable function
   $ t \colon [-\nicefrac{\ell}{2}, \nicefrac{\ell}{2}] \longrightarrow
        [-\nicefrac{\ell}{2}, \nicefrac{\ell}{2}]$, 
such that $t(-\nicefrac{\ell}{2}) = -\nicefrac{\ell}{2}$, $t(\nicefrac{\ell}{2}) =
\nicefrac{\ell}{2}$, and $t'(x) > 0$ for every $x$; then we define the perturbed doping as $\bar{C}(x) \Def \tilde{C}(t(x))$.

In order to generate $t$, we impose that $p(x) \Def t(x) - x$ is a polynomial of
degree 2 or 3. Due to the properties of the function $t$, we obtain that
$p(-\nicefrac{\ell}{2}) = 0$ and $p(\nicefrac{\ell}{2}) = 0$. We randomly decide
whether to use a polynomial of degree 2 or of degree 3, that is, we use
\[ 
    p(x) = k \left(x + \nicefrac{\ell}{2}\right) \left(x - \nicefrac{\ell}{2}\right) \qquad\text{ or } \qquad   p(x) = k \left(x + \nicefrac{\ell}{2}\right) \left(x - \nicefrac{\ell}{2}\right) (x - \alpha) 
\]
where $k$ and $\alpha$ are random constants chosen so that  
$p'(x) > -1$ for every $x$.

Applying the transformation on a doping function, we can generate new
samples for our dataset whose doping can not be described simply by choosing
some suitable parameters in~\eqref{eq:doping_space}.

\section{ResNet structure}
\label{sec:resnet-details}
Exactly as we did for the MLP models, all of our ResNets are preceeded by a down-scaling
interpolation layer and, at the end, there is an up-scaling layer that restores the
original dimension of the data; the scaling layers use cubic interpolation to describe
the signals or the dopings on different spacial grids. In this case, we fix the size
of the coarse grid to 256, i.e., to a power of two that is close to the size that we
have seen performs better in the MLP model. Using a power of two we ensure that the
downscaling blocks of the ResNet always deal with  an even number of neurons.

Following \cite{KaimingZhang2015}, the structure of our ResNets is made of three
different parts: the gate, that elaborates the input from the down-scaling layer;
the encoder, which applies the convolutional layers in order to extract the
most relevant features of the signal; and the decoder, which takes as input
the features recognized by the encoder and produces the model prediction. In the
following part, we describe in detail the structure of each part of the network
and its associated hyperparameters.

\begin{figure}[!h]
\begin{center}
\begin{tikzpicture}[
    every node/.append style={draw}
  ]
  \node[rotate=90, minimum width=2.6cm] (dscaling) at (0,0) {down-scaling};
  \node[minimum width=2cm] (gate) at (3,0) {Gate};
  \node[minimum width=2cm] (encoder) at (6,0) {Encoder};
  \node[minimum width=2cm] (decoder) at (9,0) {Decoder};
  \node[rotate=90, minimum width=2.6cm] (uscaling) at (12,0) {up-scaling};

  \path[->]
      (dscaling) edge (gate)
      (gate) edge (encoder)
      (encoder) edge (decoder)
      (decoder) edge (uscaling);

      \draw[draw=black] (1.5, -1.3) rectangle ++(9, 2.6);

  \node[draw=none] at (6, 1) {ResNet};
\end{tikzpicture}
\end{center}
\end{figure}

\paragraph{Gate} This is the first part of our network (after the down-sampling
layer). It consists of a convolutional layer, a batch normalization layer, and
an activation layer. In the convolutional layer, we choose the kernel sizes from
the set $\{3, 5, 7, 9\}$. The number of output channels, instead, is chosen from
the set $\{8, 16, 24, 32\}$. Finally, the stride of the convolution layer is
chosen from the set $\{1,2,4\}$. This leads to $4 \times 4 \times 3 = 48$
possible convolutional layers for our gate. Taking into account that the
activation layer always applies a ReLU (that does not require parameters) and
that also the normalization layer is fixed, we have a total of 54 possible
configurations for the gate.

\paragraph{Encoder} The encoder is built by stacking several blocks of the same type.
We consider two different kinds of blocks: the ``basic blocks'' described in
\cite{KaimingZhang2015} and the ``fixed channel block''. Both blocks are made of
the following layers: a convolutional layer, a batch normalization, a ReLU
activation layer, another convolution, and, finally, a normalization. The two
convolutional layers have a fixed kernel size of 3: they do not have bias
and the padding is chosen so that the size of the output is preserved
(therefore, in our case, the padding is equal to 1).

The difference between a
ResNet and a plain convolutional neural network is the fact that the input of
each block is not simply the output of the previous one. Instead, each block $B_i$
saves the input $x_i$ it receives and, after its computations, sums its output
$B_i(x_i)$ with the original input $x_i$ (or with a simple function of the input
$s_i(x_i)$). In this way, the input of the layer $B_{i + 1}$ is
\begin{equation}\label{eq:residual_sum} x_{i + 1} = x_i + B_i(x_i) . \end{equation}
Usually, in a ResNet, the shape of the data may change along the layers; indeed,
$x_i$ is a bidimensional tensor: the first index represents the spatial position
and while second one is the channel. The input of the first block of the encoder
has shape (256, $k_1$) where $k_1$ is the number of output channels of the gate but,
while the blocks become deeper, $k_i$ increases and the number of points decreases.
A block $B_i$ so that $B(x_i)$ has a different shape respect to $x_i$ is called a
\emph{downsampler} block. In our networks, each downsampler block halves the size
of the first index of the tensor; so, for example, the output of the first
downsampler block will have size (128, $k_{i}$).
The only difference between the ``basic blocks'' and the ``fixed channel blocks''
is in how the perform the downsample: a downsampling basic block increases the
number of channels by a factor two, while a fixed channel block does not.

It is worth noting that \cref{eq:residual_sum} can not be applied by the downsampling
blocks because of the different shapes of the tensors. In this case, the equation
becomes
\begin{equation*} x_{i + 1} = s_i(x_i) + B_i(x_i) . \end{equation*}
where $s_i$ is called "shortcut operation"; in our network, $s_i$ is performed
by a convolutional layer with kernel size equal to 1 and stride equal to 2
(basic blocks), or kernel size equal to 2 and stride equal to 2 (for the
fixed channel blocks), followed by a normalization block.

The convolutional layer of the shortcut of the fixed channel blocks
forbids any communication between channels, i.e., each element of the output tensor
depends only on the values of the elements of the input tensor that share the
same index for the channel (or, in other words, using PyTorch we impose that the
number of groups of the convolutional layer is equal to the number of channels).
For what concerns the computation of the output (and not the shortcut), the
dimensional reduction is obtained by setting the stride of the first
convolutional layer to 2.

Therefore, the only free parameters that we have left for our encoder are the
number of blocks, their kind, and a downsampling flag for each
block. Our encoders consist of one, two, or three
blocks of the same kind. For encoders made of basic blocks,
we allow two different configurations of the downsample flag: true for all
blocks or false for all blocks. For fixed channel blocks, we always set the
downsample flag to true. We have a total of 6 configurations for the
encoders with the basic blocks and 3 configurations that use the fixed channels
blocks, for a total of 9 possible configurations.

\begin{figure}[!h]
\begin{center}
\begin{tikzpicture}[
    every node/.append style={draw, minimum width=3cm, minimum height=0.6cm}
  ]
  \node[rotate=90] (c1) at (0, 0) {Conv. layer};
  \node[rotate=90] (b1) at (1, 0) {Batch norm.};
  \node[rotate=90] (relu) at (2, 0) {ReLU Activation};
  \node[rotate=90] (c2) at (3, 0) {Conv. layer};
  \node[rotate=90] (b2) at (4, 0) {Batch norm.};

  \draw (-1, -2) rectangle ++(6, 4);
  \node[draw=none] at (2, 2.3) {Block $i$};

  \node[minimum height=4cm, minimum width=2.2cm] (bb) at (-4, 0) {Block $i \! - \! 1$};

  \node[minimum height=4cm, minimum width=2.2cm] (ba) at (8, 0) {Block $i \! + \! 1$};

  \node[circle, minimum size=0.35] (circ1) at (-1.75, 0) {$\phantom{+}$};

  \node[circle, minimum size=0.35] (circ2) at (5.75, 0) {$+$};

  \path[->]
  (bb) edge[-] (circ1)
      (circ1) edge (c1)
      (c1) edge (b1)
      (b1) edge (relu)
      (relu) edge (c2)
      (c2) edge (b2)
      (b2) edge (circ2)
      (circ2) edge (ba)
      (circ1) edge[-] (-1.75, -3)
      (-1.75, -3) edge[-] node[pos=0.5,below,draw=none] {{\it shortcut}} (5.75, -3)
      (5.75, -3) edge (circ2);

\end{tikzpicture}
\end{center}
\end{figure}

\paragraph{Decoder} This is essentially a multilayer perceptron, and we allow 1
or 2 hidden layers. If we choose a configuration with 1 layer, the size of this
layer could be 100, 150 or 200. With 2 layers, there is a total of 15 possible
configurations obtained by choosing the size of the first layer  in $\{100,
150, 200, 250, 300\}$ and the second one in $\{100, 150, 200\}$. In total, we
have therefore 18 different configurations for the decoder.

\bibliographystyle{siamplain}        

\bibliography{references.bib}

\end{document}